\documentclass[11pt, reqno, psamsfonts]{amsart}
\pdfoutput=1

\usepackage{amssymb}
\usepackage{amsthm}
\usepackage{amsmath}
\usepackage{latexsym}
\usepackage[T1]{fontenc}
\usepackage[utf8]{inputenc}
\usepackage[russian, french, english]{babel}

\usepackage{graphicx}
\usepackage{wrapfig}
\usepackage[justification=centering, labelfont=bf]{caption}
\usepackage{mathtools}
\usepackage[hidelinks]{hyperref}
\usepackage{tikz}
\usepackage{amsbsy}
\usepackage[inline]{enumitem}
\usepackage{mathrsfs}
\usetikzlibrary{shapes,snakes}
\usetikzlibrary{arrows.meta}
\usetikzlibrary{decorations.pathreplacing}
\usepackage{float}
\usepackage{multicol}
\usepackage{stmaryrd}
\usepackage{eucal}
\usepackage{mathabx}\usepackage{lmodern}
\usepackage{upgreek}
\usepackage{titlesec}
\usepackage[spacing=true,kerning=true,babel=true,tracking=true]{microtype}
\usepackage[shortcuts]{extdash}
\usepackage[foot]{amsaddr} 
\usepackage{framed}
\usepackage{algpseudocode}
\usepackage{algorithm}
\definecolor{shadecolor}{gray}{0.95}
\usepackage[left=1in,right=1in,top=0.99in,bottom=1in,bindingoffset=0cm]{geometry}

\usepackage[backend=biber, style=alphabetic, sorting=nyt, maxnames=100]{biblatex}
\title{Equitable Colorings of Borel Graphs}
\date{}
\author{Anton~Bernshteyn}
\address[A.B.]{\normalfont School of Mathematics, Georgia Institute of Technology, Atlanta, GA, USA}
\email{bahtoh@gatech.edu}

\author{Clinton T. Conley}
\address[C.C.]{\normalfont Department of Mathematical Sciences, Carnegie Mellon University, Pittsburgh, PA, USA}
\email{clintonc@andrew.cmu.edu}

\thanks{Research of A.B.~is partially supported by the NSF grant DMS-2045412. Research of C.C.~is partially supported by the NSF Grant DMS-1855579.}

\newtheoremstyle{bfnote}%
{}{}%
{\slshape}{}%
{\bfseries}{\bfseries.}%
{ }%
{\thmname{#1}\thmnumber{ #2}\thmnote{ \ep{\normalfont{}#3}}}

\newtheoremstyle{defbfnote}%
{}{}%
{}{}%
{\bfseries}{.}%
{ }%
{\thmname{#1}\thmnumber{ #2}\thmnote{ (#3)}}

\newtheoremstyle{claim}%
{}{}%
{\slshape}{}%
{\itshape}{.}%
{ }%
{\thmname{#1}\thmnumber{ #2}\thmnote{ \ep{\normalfont{}#3}}}

\theoremstyle{bfnote}

\newtheorem{theo}[equation]{Theorem}
\newtheorem{prop}[equation]{Proposition}
\newtheorem{lemma}[equation]{Lemma}
\newtheorem{corl}[equation]{Corollary}
\newtheorem{conj}[equation]{Conjecture}

\newtheorem*{claim*}{Claim}
\newtheorem{big_claim}[equation]{Claim}
\newtheorem*{corl*}{Corollary}

\theoremstyle{claim}

\newcounter{ForClaims}[section]
\newtheorem{claim}{Claim}[ForClaims]

\newcommand*{\myproofname}{Proof}
\newenvironment{claimproof}[1][\myproofname]{\begin{proof}[#1]}{\end{proof}}

\theoremstyle{definition}
\newtheorem{defn}[equation]{Definition}
\newtheorem*{defn*}{Definition}

\newtheorem*{exmp*}{Example}

\newtheorem*{assum*}{Assumptions}

\theoremstyle{remark}
\newtheorem*{ques*}{Question}
\newtheorem*{remk*}{Remark}

\newcommand{\0}{\varnothing}
\newcommand{\set}[1]{\{#1\}}

\newcommand{\dom}{\mathrm{dom}}
\newcommand{\im}{\mathrm{im}}

\newcommand{\N}{{\mathbb{N}}}

\newcommand{\R}{\mathbb{R}}

\renewcommand{\epsilon}{\varepsilon}
\renewcommand{\phi}{\varphi}
\renewcommand{\theta}{\vartheta}
\renewcommand{\leq}{\leqslant}
\renewcommand{\geq}{\geqslant}

\newcommand{\fins}[1]{[#1]^{<\infty}}

\renewcommand{\Prob}{\mathsf{Prob}}
\newcommand{\defeq}{\coloneqq}

\newcommand{\emphd}[1]{{\fontseries{b}\selectfont\textsf{#1}}}

\newcommand{\dist}{{\mathrm{dist}}}
\newcommand{\disc}{{\mathrm{disc}}}
\newcommand{\Diff}{\mathrm{d}}

\newcommand{\recolor}{\oplus}
\newcommand{\col}[1]{\mathcal{#1}}
\newcommand{\diff}{\updelta}
\newcommand{\RM}{\mathsf{RM}}
\newcommand{\Erg}{\mathsf{Erg}}
\newcommand{\Inv}{\mathsf{Inv}}

\newcommand{\dominates}{\succcurlyeq}
\newcommand{\inj}{\preccurlyeq}
\newcommand{\surj}{\succcurlyeq}

\newcommand{\bemph}[1]{{\normalfont#1}} 
\newcommand{\ep}[1]{\bemph{(}#1\bemph{)}} 

\newenvironment{scproof}[1][]{\begin{proof}[\textsc{\upshape{Proof}}#1]}{\end{proof}}

\usepackage{xspace}

\numberwithin{equation}{section}

\renewcommand{\thesubsection}{\arabic{section}.\Alph{subsection}}

\usepackage{etoolbox}

\titleformat{\section}[block]{\scshape\filcenter}{\thesection.}{1ex}{}
\titleformat{\subsection}[block]{\bfseries\filcenter}{\thesubsection.}{1ex}{}
\titleformat{\subsubsection}[runin]{\itshape}{\bfseries\upshape\thesubsubsection.}{1ex}{}[.---]

\titlespacing*{\section}{0pt}{*3}{*1}
\titlespacing*{\subsection}{0pt}{*3}{*1}
\titlespacing*{\subsubsection}{0pt}{*1.5}{*0}

\makeatletter
\newcommand{\neutralize}[1]{\expandafter\let\csname c@#1\endcsname\count@}
\makeatother

\renewbibmacro{in:}{}

\renewbibmacro*{volume+number+eid}{%
	\printfield{volume}%
	\setunit*{\addnbspace}
	\printfield{number}%
	\setunit{\addcomma\space}%
	\printfield{eid}}

\DeclareFieldFormat[article]{volume}{\textbf{#1}\space}
\DeclareFieldFormat[article]{number}{\mkbibparens{#1}}

\DeclareFieldFormat{journaltitle}{#1,}
\DeclareFieldFormat[thesis]{title}{\mkbibemph{#1}\addperiod}
\DeclareFieldFormat[article, unpublished, thesis]{title}{\mkbibemph{#1},}
\DeclareFieldFormat[book]{title}{\mkbibemph{#1}\addperiod}
\DeclareFieldFormat[unpublished]{howpublished}{#1, }

\DeclareFieldFormat{pages}{#1}

\DeclareFieldFormat[article]{series}{Ser.~#1\addcomma}

\setlist{topsep=3pt,itemsep=3pt}

\begin{document}
	
	\pagestyle{plain}
	
	\maketitle
	
	
	
	\begin{abstract}
		Hajnal and Szemer\'edi proved that if $G$ is a finite graph with maximum degree $\Delta$, then for every integer $k \geq \Delta+1$, $G$ has a proper coloring with $k$ colors in which every two color classes differ in size at most by $1$; such colorings are called equitable. We obtain an analog of this result for infinite graphs in the Borel setting. Specifically, we show that if $G$ is an aperiodic Borel graph of finite maximum degree $\Delta$, then for each $k \geq \Delta + 1$, $G$ has a Borel proper $k$-coloring in which every two color classes are related by an element of the Borel full semigroup of $G$. In particular, such colorings are equitable with respect to every $G$-invariant probability measure. 
		We also establish a measurable version of a result of Kostochka and Nakprasit on equitable $\Delta$-colorings of graphs with small average degree. Namely, we prove that if $\Delta \geq 3$, $G$ does not contain a clique on $\Delta + 1$ vertices, and $\mu$ is an atomless $G$-invariant probability measure such that the average degree of $G$ with respect to $\mu$ is at most $\Delta/5$, then $G$ has a $\mu$-equitable $\Delta$-coloring. As steps towards the proof of this result, we establish measurable and list coloring extensions of a strengthening of Brooks's theorem due to Kostochka and Nakprasit. 
	\end{abstract}
	
	\section{Introduction}
	
	\subsection{Hajnal--Szemer\'edi theorem for Borel graphs}
	
	Let $G$ be a \ep{simple undirected} graph with vertex set $V(G)$ and edge set $E(G)$. For a vertex $x \in V(G)$, we denote the {neighborhood} of $x$ in $G$ by $N_G(x)$ and write $\deg_G(x) \defeq |N_G(x)|$ for the {degree} of $x$ in $G$. The \emphd{maximum degree} of $G$, denoted by $\Delta(G)$, is defined by $\Delta(G) \defeq \sup_{x \in V(G)} \deg_G(x)$.
	
	Given a set $\col{C}$, a \emphd{$\col{C}$-coloring} of $G$ is simply a mapping $f \colon V(G) \to \col{C}$; in this context we call the elements of $\col{C}$ \emphd{colors}. A $\col{C}$-coloring $f$ is \emphd{proper} if $f(x) \neq f(y)$ whenever $x$ and $y$ are adjacent in $G$. We will be mostly interested in the case when $\col{C}$ is finite. With a slight abuse of terminology, we refer to $\col{C}$-colorings with a fixed finite set $\col{C}$ of size $k$, say $\set{1, \ldots, k}$, as \emphd{$k$-colorings}. Here, and throughout the paper, $k$ denotes a positive integer. 
	
	Given a $\col{C}$-coloring $f$, we refer to the sets $f^{-1}(\alpha)$ for $\alpha \in \col{C}$ as the \emphd{color classes} of $f$. A proper $\col{C}$-coloring of a {finite} graph $G$ is \emphd{equitable} if every two color classes of $f$ differ in size at most by $1$. In particular, if $|V(G)|$ is divisible by $k \geq 1$, then in an equitable $k$-coloring of $G$ all color classes must be of size precisely $|V(G)|/k$. In contrast to ordinary coloring, a graph with an equitable $k$-coloring need not have an equitable $(k+1)$-coloring. Nevertheless, Erd\H{o}s \cite{ErdConj} conjectured that every finite graph of maximum degree $\Delta$ has an equitable $k$-coloring for each $k \geq \Delta+1$. Erd\H{o}s's conjecture was confirmed by Hajnal and Szemer\'edi~\cite{HSz}:
	
	\begin{theo}[{Hajnal--Szemer\'edi \cite{HSz}}]\label{theo:HSz}
		Let $G$ be a finite graph of maximum degree $\Delta$. If $k \geq \Delta+1$, then $G$ has an equitable $k$-coloring.
	\end{theo}
	
	The original proof of Theorem~\ref{theo:HSz} due to Hajnal and Szemer\'edi was surprisingly difficult, but it was significantly simplified in the work of Mydlarz and Szemer\'edi \ep{unpublished, see \cite{equit}} and Kierstead and Kostochka \cite{KierKos}, culminating in a two-page proof. Moreover, their argument provides an efficient algorithm that builds a desired equitable coloring \cite{equit}.
	
	The central result of this paper is an extension of Theorem~\ref{theo:HSz} to equitable colorings of {infinite} graphs. 
	Specifically, if $G$ is a graph whose vertex set $V(G)$ carries a probability measure, then it is natural to call a proper $k$-coloring $f$ of $G$ {equitable} if every color class of $f$ has measure $1/k$. Notice that, in order for this definition to be sensible, we must require that every color class of $f$ is a {measurable subset} of $V(G)$. Questions regarding the behavior of colorings, matchings, and other combinatorial constructions under extra measurability constraints are studied in the area of \emph{descriptive combinatorics}, which has attracted considerable attention in recent years; see \cite{KechrisMarks} for a comprehensive survey.
	
	Before stating our results, we need to introduce some relevant terminology. Our main references for descriptive set theory are \cite{KechrisDST, AnushDST}. By a \emphd{Borel graph} we mean a graph $G$ whose vertex set $V(G)$ is a standard Borel space and whose edge set $E(G)$ is a Borel subset of $V(G) \times V(G)$. 
	If $G$ is a Borel graph and $\col{C}$ is a standard Borel space, then a $\col{C}$-coloring $f \colon V(G) \to \col{C}$ is \emphd{Borel} if it is a Borel function, i.e., if preimages of Borel subsets of $\col{C}$ under $f$ are Borel in $V(G)$. When $\col{C}$ is countable, this is equivalent to saying that every color class of $f$ is a Borel subset of $V(G)$. The smallest cardinality of a standard Borel space $\col{C}$ such that $G$ admits a Borel proper $\col{C}$-coloring is called the \emphd{Borel chromatic number} of $G$ and is denoted by $\chi_{\mathrm{B}}(G)$. Similarly, given a probability measure $\mu$ on $V(G)$, we can talk about \emphd{$\mu$-measurable} colorings $f \colon V(G) \to \col{C}$ \ep{i.e., such that $f$-preimages of Borel subsets of $\col{C}$ are $\mu$-measurable} and define the \emphd{$\mu$-measurable chromatic number} $\chi_\mu(G)$ of $G$ as the smallest cardinality of a standard Borel space $\col{C}$ such that $G$ admits a $\mu$-measurable proper $\col{C}$-coloring. Borel chromatic numbers were first introduced and systematically studied by Kechris, Solecki, and Todorcevic in their seminal paper \cite{KST}. Among several other results, they established the following fact:
	
	\begin{theo}[{Kechris--Solecki--Todorcevic \cite[Proposition 4.6]{KST}}]\label{theo:KST}
		Let $G$ be a Borel graph of finite maximum degree $\Delta$. Then $\chi_{\mathrm{B}}(G) \leq \Delta+1$.
	\end{theo}

	In view of Theorem~\ref{theo:KST}, it is meaningful to ask for Borel $(\Delta+1)$-colorings with extra properties \ep{such as being equitable}. We remark that, according to a startling result of Marks~\cite{Marks}, the upper bound $\chi_{\mathrm{B}}(G) \leq \Delta + 1$ is sharp, even for {acyclic} graphs $G$.
	
	\begin{defn}\label{defn:mu-equit}
		Let $G$ be a Borel graph and let $\mu$ be a probability measure on $V(G)$. A~\emphd{$\mu$\=/equitable $k$-coloring} of $G$ is a $\mu$-measurable proper $k$-coloring $f$ of $G$ such that $\mu(f^{-1}(\alpha)) = 1/k$ for every color $\alpha$.
	\end{defn}

	Just as the definition of equitable coloring for finite graphs uses the \ep{normalized} counting measure on $V(G)$, Definition~\ref{defn:mu-equit} is most natural for measures $\mu$ that ``assign the same weight'' to every vertex of $G$. Formally, let $\llbracket G \rrbracket$ denote the \emphd{Borel full semigroup}
	of $G$, i.e., the set of all Borel bijections $\phi \colon A \to B$, where $A$ and $B$ are Borel subsets of $V(G)$, such that for all $x \in A$, $\phi(x)$ and $x$ are joined by a path in $G$. We say that Borel subsets $A$, $B \subseteq V(G)$ are \emphd{$G$\=/equidecomposable}, in symbols $A \approx_G B$, if there is $\phi \in \llbracket G \rrbracket$ with $\dom(\phi) = A$ and $\im(\phi) = B$. A probability measure $\mu$ on $V(G)$ is said to be \emphd{$G$-invariant} if $\mu(A) = \mu(B)$ whenever $A \approx_G B$.  When the maximum degree of $G$ is finite, this is equivalent to the following ``double\-/counting'' identity:
	\[
	\int_A |N_G(x) \cap B| \,\Diff\mu(x) \,=\, \int_B |N_G(y) \cap A| \,\Diff \mu(y), \qquad \text{for all Borel } A,\, B \subseteq V(G). 
	\]
	
	\begin{defn}
		Let $G$ be a Borel graph. A \emphd{Borel\-/equitable $k$-coloring} of $G$ is a Borel proper $k$-coloring $f$ of $G$ such that $f^{-1}(\alpha) \,\approx_G\, f^{-1}(\beta)$ for every pair of colors $\alpha$ and $\beta$.
	\end{defn}

	It follows immediately that a Borel\-/equitable $k$-coloring of $G$ is $\mu$-equitable with respect to {every} $G$-invariant probability measure $\mu$. 
	
	We prove a version of the \hyperref[theo:HSz]{Hajnal--Szemer\'edi theorem} for Borel-equitable colorings. In order to avoid divisibility issues and thus make its statement simpler, we shall focus on \emphd{aperiodic} graphs, i.e., those in which every connected component is infinite. However, the extension to graphs with finite components is straightforward; see, e.g., Lemma~\ref{lemma:meas_HSz}.
	
	\begin{theo}[\textls{Borel Hajnal--Szemer\'edi}]\label{theo:BHSz}
		Let $G$ be an aperiodic Borel graph of finite maximum degree $\Delta$. If  $k \geq \Delta + 1$, then $G$ has a Borel\-/equitable $k$-coloring.
	\end{theo}

	Additionally, we show that if a given coloring $f$ is ``approximately equitable,'' then $f$ is actually ``close'' to an equitable coloring. To make this precise, we need a couple more definitions. A subset $A \subseteq V(G)$ is \emphd{$G$-invariant} if it is a union of connected components of $G$, i.e., if no edge of $G$ joins a vertex in $A$ to a vertex in $V(G) \setminus A$. A $G$-invariant measure $\mu$ is \emphd{$G$-ergodic} \ep{or simply \emphd{ergodic} if $G$ is clear from the context} if every $G$-invariant Borel subset of $V(G)$ is either $\mu$-null or $\mu$-conull. Every $G$-invariant probability measure can be decomposed as a convex combination of ergodic measures; see Theorem~\ref{theo:UED} for a precise statement of this fact. Now, fix a nonempty finite color set $\col{C}$ and let $\mu$ be a probability measure on $V(G)$. The \emphd{$\mu$-discrepancy} of a $\mu$-measurable $\col{C}$-coloring $f$, in symbols $\disc_\mu(f)$, is defined by the formula
	\[
	\disc_\mu(f) \,\defeq\, \max_{\alpha \in \col{C}} \left|\mu(f^{-1}(\alpha)) \,-\, |\col{C}|^{-1}\right|.
	\]
	The \emphd{$\mu$-distance} between two $\col{C}$-colorings $f$, $g$, in symbols $\dist_\mu(f,g)$, is defined as
	\[
	\dist_\mu(f, g) \,\defeq\, \mu\left(\set{x \in V(G) \,:\, f(x) \neq g(x)}\right).
	\]
	We establish the following strengthening of Theorem~\ref{theo:BHSz}:
	
	\begin{theo}[\textls{Stable Borel Hajnal--Szemer\'edi}]\label{theo:stable_BHSz}
		Let $G$ be an aperiodic Borel graph of finite maximum degree $\Delta$ and let $f$ be a Borel proper $k$-coloring of $G$, where $k \geq \Delta + 1$.
		\begin{enumerate}[label=\ep{\itshape\alph*}]
		    \item\label{item:stable_a} For every $G$-invariant probability measure $\mu$, there is a $\mu$\-/equitable $k$-coloring $g$ such that
		    \begin{equation}\label{eq:stable}
			\dist_\mu(f, g) \,\leq\, 7^{k+1} \cdot \disc_\mu(f).
		    \end{equation}
		    \item\label{item:stable_b} Furthermore, $G$ has a Borel\-/equitable $k$-coloring $g$ such that \eqref{eq:stable} holds for every ergodic $G$-invariant probability measure $\mu$ simultaneously.
		\end{enumerate}
	\end{theo}
	
	We did not make an effort to optimize the coefficient $7^{k+1}$ in front of $\disc_\mu(f)$ in \eqref{eq:stable}; however, even with more care our proof techniques seem to lead to exponential dependence on $k$. It would be interesting to know whether exponential dependence is necessary; in principle, it might even be possible to replace $7^{k+1}$ by a {linear} function of $k$.
	

	\subsection{Equitable $\Delta$-colorings}
	
	By Brooks's theorem \cite[Theorem~5.2.4]{Die00}, ``most'' connected finite graphs of maximum degree $\Delta$ can be properly colored using only $\Delta$ colors; the only exceptions are odd cycles and complete graphs.
	For equitable $\Delta$-colorings, at least one new class of pathological examples is known: the complete bipartite graphs $K_{\Delta,\Delta}$ for odd $\Delta$. 
	The following analog of Brooks's theorem was conjectured by Chen, Lih, and Wu \cite{CLW}: 
	
	\begin{conj}[{Chen--Lih--Wu \cite{CLW}}]\label{conj:equit_Brooks}
		Let $G$ be a connected finite graph of maximum degree $\Delta \geq 1$. Then $G$ has an equitable $\Delta$-coloring, unless:
		\begin{quote}
			\begin{enumerate*}[label=\ep{\normalfont{}\itshape\alph*},itemjoin={\quad}]
				\item $\Delta=2$ and $G$ is an odd cycle; 
				
				\item $G \cong K_{\Delta+1}$;\quad or 
				
				\item $\Delta$ is odd and $G \cong K_{\Delta,\Delta}$.
			\end{enumerate*}
		\end{quote}
	\end{conj}
	
	To date, Conjecture~\ref{conj:equit_Brooks} remains open. However, some partial results are known; see \cite{Lih_survey} for a survey. In particular, Kostochka and Nakprasit \cite{KN} proved that $G$ has an equitable $\Delta$-coloring provided that the \emphd{average degree} of $G$, i.e., the quantity \[d(G) \,\defeq\, \frac{\sum_{x \in V(G)} \deg_G(x)}{|V(G)|} \,=\, \frac{2|E(G)|}{|V(G)|},\]  is considerably smaller than $\Delta$:
	
	\begin{theo}[{Kostochka--Nakprasit \cite[Theorem 1]{KN}}]\label{theo:KN}
		Let $G$ be a finite graph of maximum degree $\Delta \geq 46$ and without a clique on $\Delta+1$ vertices. If the average degree of $G$ is at most $\Delta/5$, then $G$ admits an equitable $\Delta$-coloring.
	\end{theo}

	Our second main result is a measurable version of Theorem~\ref{theo:KN}. Unfortunately, Brooks's theorem {fails} in the setting of Borel colorings \ep{as mentioned earlier, Marks \cite{Marks} showed that the bound $\chi_{\mathrm{B}}(G) \leq \Delta + 1$ is sharp even for acyclic graphs $G$}. However, Conley, Marks, and Tucker-Drob established a version of Brooks's theorem for \emph{measurable} colorings:
	
	\begin{theo}[\textls{Measurable Brooks}; {Conley--Marks--Tucker-Drob \cite[Theorem 1.2(1)]{CMTD}}]\label{theo:meas_Brooks}
		Let $G$ be a Borel graph of finite maximum degree $\Delta \geq 3$ and without a clique on $\Delta + 1$ vertices. Then $\chi_\mu(G) \leq \Delta$ for every probability measure $\mu$ on $V(G)$.
	\end{theo}
	
	Thus, there is still hope of constructing $\mu$-equitable $\Delta$-colorings. For a Borel graph $G$ of finite maximum degree and a probability measure $\mu$ on $V(G)$, let the \emphd{$\mu$-average degree} of $G$ be
	\[
	d_\mu(G) \,\defeq\, \int_{V(G)} \deg_G(x) \,\Diff \mu(x).
	\]
	The reader familiar with measurable graph theory would notice that $d_\mu(G) = 2\mathsf{C}_\mu(G)$, where $\mathsf{C}_\mu(G)$ is the \emph{cost} of $G$ \ep{see \cite[Chapter~III]{KechrisMiller}}. Recall that a measure $\mu$ is called \emphd{atomless} if $\mu(\set{x}) = 0$ for every point $x$. \ep{In particular, if $G$ is an aperiodic Borel graph, then every $G$-invariant probability measure on $V(G)$ is atomless.} We prove the following measurable analog of Theorem~\ref{theo:KN}:
	
	\begin{theo}\label{theo:meas_KN}
		Let $G$ be a Borel graph of finite maximum degree $\Delta \geq 3$ and without a clique on $\Delta + 1$ vertices. If $\mu$ is an atomless $G$-invariant probability measure on $V(G)$ such that $d_\mu(G) \leq \Delta/5$, then $G$ has a $\mu$-equitable $\Delta$-coloring.
	\end{theo}

	It is possible that the upper bound on $d_\mu(G)$ in Theorem~\ref{theo:meas_KN} is not actually necessary. Indeed, we suspect that the following version of Conjecture~\ref{conj:equit_Brooks} should hold in full generality:
	
	\begin{conj}
		Let $G$ be an aperiodic Borel graph of finite maximum degree $\Delta \geq 3$ and let $\mu$ be a $G$-invariant probability measure on $V(G)$. Then $G$ has a $\mu$-equitable $\Delta$-coloring.
	\end{conj}

	\subsection{Domination for partial $\Delta$-colorings}
	
	In order to prove Theorem~\ref{theo:KN}, Kostochka and Nakprasit established a useful auxiliary result concerning the relationship between $\Delta$-colorings of a finite graph $G$ and those of its subgraphs. 
	Let $G$ be a graph and let $\col{C}$ be a set. A \emphd{partial $\col{C}$-coloring} of $G$ is a function $f \colon U \to \col{C}$ with $U \subseteq V(G)$; to indicate that $f$ is a partial $\col{C}$-coloring, we write $f \colon V(G) \rightharpoonup \col{C}$. 
	A partial $\col{C}$-coloring $f$ is \emphd{proper} if $f(x) \neq f(y)$ whenever $x$, $y \in \dom(f)$ are adjacent. Given partial $\col{C}$-colorings $f$, $g$ of a {finite} graph $G$, we say that $f$ \emphd{dominates} $g$, in symbols $f \dominates g$, if $|f^{-1}(\alpha)| \geq |g^{-1}(\alpha)|$ for all $\alpha \in \col{C}$. 
	In particular, if $f$ is an {extension} of $g$ \ep{i.e., if $f \supseteq g$}, then $f$ dominates $g$; but in general the relation $f \dominates g$ says nothing about the values $f$ and $g$ take at individual vertices.
	
	\begin{theo}[{Kostochka--Nakprasit \cite[Theorem 2]{KN}}]\label{theo:KN_domination}
		Let $G$ be a finite graph of maximum degree $\Delta \geq 3$ and without a clique on $\Delta + 1$ vertices. If $g$ is a proper partial $\Delta$-coloring of $G$, then $G$ has a proper $\Delta$-coloring $f$ such that $f \dominates g$.
	\end{theo}

	Unsurprisingly, our proof of Theorem~\ref{theo:meas_KN} similarly relies on a measurable version of Theorem~\ref{theo:KN_domination}. To state it, we extend the notion of domination for partial colorings to the measurable context in the obvious way. Namely, if $G$ is a Borel graph and $\mu$ is a probability measure on $V(G)$, then, given a pair of $\mu$-measurable partial colorings $f$, $g$, we say that $f$ \emphd{dominates} $g$ \emphd{with respect to} $\mu$, in symbols $f \dominates_\mu g$, if $\mu(f^{-1}(\alpha)) \geq \mu(g^{-1}(\alpha))$ for every color $\alpha$.
	
	\begin{theo}[\textls{Measurable domination}]\label{theo:meas_domination}
		Let $G$ be a Borel graph of finite maximum degree $\Delta \geq 3$ and without a clique on $\Delta + 1$ vertices and let $\mu$ be a $G$-invariant probability measure on $V(G)$. If $g$ is a $\mu$-measurable proper partial $\Delta$-coloring of $G$, then $G$ has a $\mu$-measurable proper $\Delta$-coloring $f$ such that $f \dominates_\mu g$.
	\end{theo}
	
	In combination with Theorem~\ref{theo:BHSz}, Theorem~\ref{theo:meas_domination} yields the following corollary:
	
	\begin{corl}[\textls{Almost equitable $\Delta$-colorings}]\label{corl:almost_equit}
		Let $G$ be an aperiodic Borel graph of finite maximum degree $\Delta \geq 3$ and let $\mu$ be a $G$-invariant probability measure on $V(G)$. Then $G$ has a $\mu$-measurable proper $\Delta$-coloring $f$ such that $\mu(f^{-1}(\alpha)) \geq 1/(\Delta + 1)$ for every color $\alpha$.
	\end{corl}
	\begin{scproof}
		By Theorem~\ref{theo:BHSz}, $G$ has a $\mu$-equitable $(\Delta + 1)$-coloring $h$. Fix an arbitrary color $\beta$ and let $g$ be the partial $\Delta$-coloring of $G$ obtained from $h$ by uncoloring all the vertices in $h^{-1}(\beta)$. Then every color class of $g$ has measure precisely $1/(\Delta + 1)$, and thus applying Theorem~\ref{theo:meas_domination} to this $g$ yields the desired $\Delta$-coloring $f$.
	\end{scproof}

	It turns out that, in order to prove Theorem~\ref{theo:meas_domination}, we must first strengthen Theorem~\ref{theo:KN_domination} by extending it to the \emph{list coloring} context. This phenomenon is not uncommon in graph coloring theory; for example, the proof of Theorem~\ref{theo:meas_Brooks} due to Conley, Marks, and Tucker-Drob relies in a similar fashion on the list coloring version of Brooks's theorem \ep{see Theorem~\ref{theo:list_Brooks} below}. As we believe our list coloring analog of Theorem~\ref{theo:KN_domination} to be of independent interest, we describe it here.
	
	List coloring is a generalization of graph coloring that was introduced independently by Vizing \cite{Vizing} and Erd\H{o}s, Rubin, and Taylor~\cite{ERT}. A \emphd{list assignment} for a graph $G$ is a mapping $\col{L}$ that assigns to each vertex $x \in V(G)$ a set $\col{L}(x)$, called the \emphd{list} of $x$. An \emphd{$\col{L}$-coloring} of $G$ is a function $f$ with domain $V(G)$ such that $f(x) \in \col{L}(x)$ for all $x \in V(G)$; similarly, a \emphd{partial $\col{L}$-coloring} is a function $f$ with $\dom(f) \subseteq V(G)$ such that $f(x) \in \col{L}(x)$ for all $x \in \dom(f)$. A \ep{partial} $\col{L}$-coloring $f$ is \emphd{proper} if $f(x) \neq f(y)$ whenever $x$ and $y$ are adjacent. Note that ordinary graph coloring is a special case of list coloring with all lists being the same.
	
\begin{figure}[h]
		\centering
		\begin{tikzpicture}
		\node[draw,regular polygon,regular polygon sides=5,minimum size=2cm] (a) at (0,0) {};
		\draw (a.corner 1) -- (a.corner 3) -- (a.corner 5) -- (a.corner 2) -- (a.corner 4) -- (a.corner 1);
		
		\node[draw,regular polygon,regular polygon sides=7,minimum size=1.5cm] (b) at (1.27,1) {};
		
		\node[draw,regular polygon,regular polygon sides=3,minimum size=1cm] (c) at (0.6,-1.3) {};
		
		\node[draw,regular polygon,regular polygon sides=4,minimum size=1.5cm] (d) at (2.13,-0.21) {};
		\draw (d.corner 1) -- (d.corner 3) (d.corner 2) -- (d.corner 4);
		
		\draw (a.corner 2) -- (-1.4,1.3);
		\draw (a.corner 2) -- (-1.6,-0.3);
		
		\draw (d.corner 1) -- (3,1);
		\draw (3,1) -- (3.5,1.5);
		\draw (3.5,1.5) -- (4,1.8);
		\draw (4,1.8) -- (4.5, 2);
		\draw[dashed] (4.5,2) -- (5.5,2.2);
		\draw (3.5,1.5) -- (4,1.2);
		\draw (4,1.2) -- (4.5, 1);
		\draw[dashed] (4.5,1) -- (5.5,0.8);
		\draw (4,1.8) -- (4.5, 1.6);
		\draw[dashed] (4.5,1.6) -- (5.5,1.4);
		
		\draw (3,1) -- (3.5,0.5);
		\draw (3.5,0.5) -- (4,0.2);
		\draw (4,0.2) -- (4.5, 0);
		\draw[dashed] (4.5,0) -- (5.5,-0.2);
		\filldraw (3,1) circle (2pt);
		\filldraw (3.5,1.5) circle (2pt);
		\filldraw (4,1.8) circle (2pt);
		\filldraw (4.5,2) circle (2pt);
		\filldraw (3,1) circle (2pt);
		\filldraw (3.5,0.5) circle (2pt);
		\filldraw (4,0.2) circle (2pt);
		\filldraw (4.5,0) circle (2pt);
		\filldraw (4,1.2) circle (2pt);
		\filldraw (4.5,1) circle (2pt);
		\filldraw (4.5,1.6) circle (2pt);
		
		\filldraw (a.corner 1) circle (2pt);
		\filldraw (a.corner 2) circle (2pt);
		\filldraw (a.corner 3) circle (2pt);
		\filldraw (a.corner 4) circle (2pt);
		\filldraw (a.corner 5) circle (2pt);
		
		\filldraw (b.corner 1) circle (2pt);
		\filldraw (b.corner 2) circle (2pt);
		\filldraw (b.corner 3) circle (2pt);
		\filldraw (b.corner 5) circle (2pt);
		\filldraw (b.corner 6) circle (2pt);
		\filldraw (b.corner 7) circle (2pt);
		
		\filldraw (c.corner 2) circle (2pt);
		\filldraw (c.corner 3) circle (2pt);
		
		\filldraw (d.corner 1) circle (2pt);
		\filldraw (d.corner 3) circle (2pt);
		\filldraw (d.corner 4) circle (2pt);
		
		\filldraw (-1.4,1.3) circle (2pt);
		\filldraw (-1.6,-0.3) circle (2pt);
		\end{tikzpicture}
		\caption{A fragment of an infinite Gallai tree.}\label{fig:Gallai}
	\end{figure}
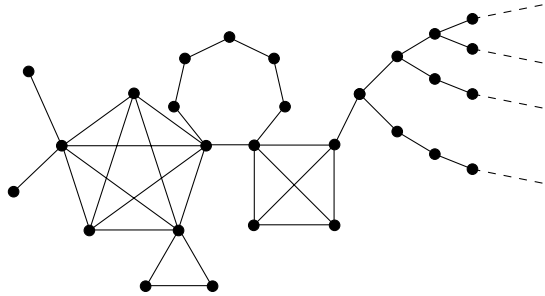

	A list assignment $\col{L}$ for a graph $G$ is called a \emphd{degree-list assignment} if $|\col{L}(x)| \geq \deg_G(x)$ for all $x \in V(G)$. A fundamental result of Borodin~\cite{Borodin} and Erd\H{o}s, Rubin, and Taylor~\cite{ERT}, which can be seen as an extension of Brooks's theorem to list coloring, provides a complete characterization of graphs $G$ that are not $\col{L}$-colorable with respect to some degree-list assignment $\col{L}$. To state it, we need to recall a few definitions. Given a vertex $x \in V(G)$, we write $G-x$ for the subgraph of $G$ obtained by deleting $x$ and all the edges incident to $x$. A \emphd{cut-vertex} in a connected graph $G$ is a vertex $x \in V(G)$ such that $G - x$ has at least two connected components. A \emphd{block} in a graph $G$ is a maximal subgraph of $G$ without a cut-vertex. A connected graph $G$ is called a \emphd{Gallai tree} if every block in $G$ is a clique or an odd cycle \ep{see Figure~\ref{fig:Gallai}}. 
	
	\begin{theo}[\textls{Brooks for list coloring}; {Borodin~\cite{Borodin}; Erd\H{o}s--Rubin--Taylor~\cite{ERT}}]\label{theo:list_Brooks}
		Let $G$ be a connected finite graph that is not a Gallai tree and let $\col{L}$ be a degree\-/list assignment for $G$. Then $G$ has a proper $\col{L}$-coloring.
	\end{theo}
	
	If $f$ and $g$ are partial $\col{L}$-colorings of a finite graph $G$, we say that $f$ \emphd{dominates} $g$, in symbols $f \dominates g$, if $|f^{-1}(\alpha)| \geq |g^{-1}(\alpha)|$ for all $\alpha \in \bigcup \set{\col{L}(x) \,:\, x \in V(G)}$. At first glance, this notion of domination may seem too strong, since each color $\alpha$ is only available to a subset of the vertices. Nevertheless, it turns out that this is the right notion for strengthening Theorem~\ref{theo:list_Brooks} and extending Theorem~\ref{theo:KN_domination} to the list coloring framework:
	
	\begin{theo}[\textls{Domination for list coloring}]\label{theo:list_domination}
		Let $G$ be a connected finite graph that is not a Gallai tree and let $\col{L}$ be a degree\-/list assignment for $G$. Suppose that $g$ is a partial proper $\col{L}$-coloring of $G$. Then $G$ has a proper $\col{L}$-coloring $f$ with $f \dominates g$.
	\end{theo}
	
	\subsection{Outline of the remainder of the paper}
	
	After the preliminary Section~\ref{sec:prelim}, we proceed to prove Theorem~\ref{theo:stable_BHSz} \ep{and thus also Theorem~\ref{theo:BHSz}} in Section~\ref{sec:BHSz}. The bulk of Section~\ref{sec:BHSz} is concerned with building {$\mu$-equitable} colorings for a fixed probability measure $\mu$; the tools that are then used to obtain a complete proof of Theorem~\ref{theo:stable_BHSz} in the purely Borel setting are the \emph{Uniform Ergodic Decomposition Theorem} of Farrell and Varadarajan \ep{see Theorem~\ref{theo:UED}}, which helps us treat all $G$-invariant probability measures simultaneously, and the properties of \emph{compressible graphs} \ep{see \S\ref{subsec:compressible}}, which allow us to handle the case when there are no $G$-invariant probability measures to begin with. In Section~\ref{sec:domination} we prove Theorem~\ref{theo:list_domination} and then use it to deduce Theorem~\ref{theo:meas_domination}. A key role in the derivation of Theorem~\ref{theo:meas_domination} is played by the method of \emph{one-ended subforests} developed by Conley, Marks, and Tucker-Drob in \cite{CMTD}. Finally, in Section~\ref{sec:meas_KN}, we prove Theorem~\ref{theo:meas_KN}.
	
	\subsubsection*{Acknowledgments}
	
	We are very grateful to Ruiyuan \ep{Ronnie} Chen for insightful discussions and, in particular, for drawing our attention to some of the facts concerning compressible graphs that are mentioned in \S\ref{subsec:compressible}. We are also grateful to the anonymous referees for their helpful suggestions.
	
	\section{Preliminaries}\label{sec:prelim}
	
	
	
	\subsubsection*{Finite sets}\label{subsec:finite}
	
	For a set $A$, let $\fins{A}$ denote the set of all finite subsets of $A$. If $X$ is a standard Borel space, then $\fins{X}$ also carries a natural standard Borel structure. One way to see this is to fix a Borel linear ordering on $X$ \ep{such an ordering exists since, by \cite[Theorem 15.6]{KechrisDST}, $X$ is isomorphic to a Borel subset of $\R$} and identify $\fins{X}$ with the set of all strictly increasing finite sequences of elements of $X$. It is a useful observation that there exists a Borel map $p \colon \fins{X}\setminus \set{\0} \to X$ such that $p(S) \in S$ for all $S \in \fins{X} \setminus \set{\0}$; for example, for a fixed Borel linear ordering on $X$, the function $S \mapsto \min S$ works.
	
	\subsubsection*{Graphs}
	
	We say that a graph $G$ is \emphd{locally countable} \ep{resp.\ \emphd{locally finite}} if every vertex of $G$ has countably \ep{resp.\ finitely} many neighbors. For $U \subseteq V(G)$, $N_G(U)$ denotes the \emphd{neighborhood} of $U$ in $G$, i.e., the set of all vertices of $G$ that have a neighbor in $U$, and $G[U]$ denotes the \emphd{subgraph of $G$ induced by $U$}, i.e., the graph with vertex set $U$ and edge set $\set{(x,y) \in E(G) \,:\, x,\, y \in U}$. A \emphd{connected component} of $G$ is a maximal subset $C \subseteq V(G)$ such that $G[C]$ is a connected graph.

	We shall use the following standard extension of Theorem~\ref{theo:KST}:
	
	\begin{prop}[{ess.~Kechris--Solecki--Todorcevic \cite[Proposition 4.6]{KST}}]\label{prop:list}
	Let $G$ be a locally countable Borel graph such that $\chi_{\mathrm{B}}(G) \leq \aleph_0$ and let $\col{C}$ be a standard Borel space of colors. Let $\col{L} \colon V(G) \to \fins{\col{C}}$ be a Borel list assignment. If $g$ is a Borel proper partial $\col{L}$-coloring of $G$, then $G$ has a Borel inclusion-maximal proper partial $\col{L}$-coloring $f$ such that $f \supseteq g$.
	\end{prop}
	\begin{scproof}
		We may assume that $g = \0$, otherwise passing to the subgraph $G[V(G) \setminus \dom(g)]$ and replacing $\col{L}$ by the list assignment $x \mapsto \col{L}(x) \setminus \set{g(y) \,:\, y \in N_G(x) \cap \dom(g)}$. Let $\theta \colon V(G) \to \N$ be a Borel proper coloring. Fix a Borel linear ordering on $\col{C}$ and recursively define partial $\col{L}$-colorings $f_n \colon \theta^{-1}(n) \rightharpoonup \col{C}$ as follows: Set 
		\[
			\col{L}_n(x) \,\defeq\, \col{L}(x) \setminus \set{\alpha \,:\, \text{there is } y \in N_G(x) \text{ with } \theta(y) < n \text{ and } f_{\theta(y)}(y) = \alpha},
		\]
		and let $f_n(x)$ be the smallest color in $\col{L}_n(x)$ if $\col{L}_n(x) \neq \0$, leaving $f_n(x)$ undefined otherwise. Then the union $f \defeq \bigcup_{n=0}^\infty f_n$ is a Borel inclusion-maximal proper partial $\col{L}$-coloring of $G$, as desired.
	\end{scproof}

	It is a useful observation that if $f$ is an inclusion-maximal proper partial $\col{C}$-coloring of a graph $G$, then each vertex $x \in V(G) \setminus \dom(f)$ has at least one neighbor of every color $\alpha \in \col{C}$ and, in particular, $\deg_G(x) \geq |\col{C}|$. Combining this observation with Proposition~\ref{prop:list}, we obtain the following:

	\begin{corl}\label{corl:extension}
		Let $G$ be a Borel graph of finite maximum degree $\Delta$ and let $k \geq \Delta + 1$. If $g$ is a Borel proper partial $k$-coloring of $G$, then $G$ has a Borel proper $k$-coloring $f$ such that $f \supseteq g$.
	\end{corl}
	\begin{scproof}
		Immediate from Proposition~\ref{prop:list} \ep{note that $\chi_\mathrm{B}(G) \leq \Delta + 1 < \aleph_0$ by Theorem~\ref{theo:KST}}. 
	\end{scproof}

	A subset $I \subseteq V(G)$ is \emphd{$G$-independent} if $I \cap N_G(I) = \0$ \ep{thus, every color class in a proper coloring of $G$ is $G$-independent}. The following is a variation of \cite[Proposition 4.2]{KST}:
	
	\begin{corl}[{ess.~Kechris--Solecki--Todorcevic \cite[Proposition 4.2]{KST}}]\label{corl:max_indep}
		Let $G$ be a locally countable Borel graph such that $\chi_{\mathrm{B}}(G) \leq \aleph_0$ and let $J \subseteq V(G)$ be a Borel $G$-independent set. Then there is a Borel maximal $G$-independent set $I \subseteq V(G)$ with $I \supseteq J$.
	\end{corl}
	\begin{scproof}
		Apply Proposition~\ref{prop:list} with $\col{L} \colon x \mapsto \set{0}$ and $g \colon J \to \set{0} \colon x \mapsto 0$.
	\end{scproof}


	\subsubsection*{Measures}
	
	Let $X$ be a standard Borel space. 
	We use $\Prob(X)$ to denote the set of all probability measures on $X$. We equip $\Prob(X)$ with the $\sigma$-algebra generated by the maps $\mu \mapsto \mu(A)$, where $A$ is a Borel subset of $X$; this makes $\Prob(X)$ into a standard Borel space \cite[\S{}17.E]{KechrisDST}.
	
	Let $G$ be a locally countable Borel graph. We write $\Inv(G)$ \ep{resp.~$\Erg(G)$} to denote the set of all $G$-invariant \ep{resp.~$G$-invariant and ergodic} probability measures on $V(G)$. Note that $\Inv(G)$ is a Borel subset of $\Prob(V(G))$, while $\Erg(G)$ is a Borel subset of $\Inv(G)$ \cite[Theorem 4.10]{KechrisCBER}. A function on $V(G)$ is \emphd{$G$-invariant} if it is constant on each connected component of $G$. Hence, a subset $A \subseteq V(G)$ is $G$-invariant if and only if its characteristic function is $G$-invariant. Observe that if $\mu \in \Inv(G)$ and $A \subseteq V(G)$ is $\mu$-null, then $A \subseteq B$ for some $G$-invariant $\mu$-null Borel set $B$. 
	The following result plays a crucial role in our proof of Theorem~\ref{theo:BHSz}:
	
	\begin{theo}[{\textls{Uniform Ergodic Decomposition}; Farrell \cite{Far}, Varadarajan \cite{Var}; see \cite[Theorem~4.11]{KechrisCBER}}]\label{theo:UED}
		Let $G$ be a locally countable Borel graph. If $\Inv(G) \neq \0$, then $\Erg(G) \neq \0$ and there is a surjective $G$-invariant Borel mapping $V(G) \to \Erg(G) \colon x \mapsto \mu_x$ such that:
		\begin{enumerate}[label=\ep{\normalfont{}E\arabic*}]
			\item for each $\mu \in \Erg(G)$, the set $\set{x \in V(G) \,:\, \mu_x = \mu}$ is $\mu$-conull;
			
			\item\label{item:E2} if $\mu \in \Inv(G)$, then $\mu = \int_X \mu_x \,\Diff \mu(x)$.
		\end{enumerate}
	\end{theo}

	For more information about this result, and about invariant measures in general, see \cite[\S{}4]{KechrisCBER}.
	
	\section{Proof of the Hajnal--Szemer\'edi theorem for Borel graphs}\label{sec:BHSz}
	
	\subsection{Recoloring moves and perfect colorings}\label{subsec:move}
	
	\begin{assum*}[for \S\ref{subsec:move}]
		 Fix a Borel graph $G$ of finite maximum degree $\Delta$ with vertex set $V$ and edge set $E$; a finite set of colors $\col{C}$ of size $k \geq \Delta + 1$; and $\mu \in \Inv(G)$.
	\end{assum*}

	A \emphd{recoloring move} is a  partial map $\phi \colon V \rightharpoonup \col{C}$ such that $\dom(\phi)$ is a nonempty finite set contained in a single connected component of $G$. The set of all recoloring moves is denoted by $\RM$. If we view each recoloring move as a finite subset of $V \times \col{C}$, then $\RM$ becomes a Borel subset of $\fins{V \times \col{C}}$. Given $m \in \N^+$, let $\RM_m$ denote the \ep{Borel} set of all recoloring moves $\phi$ with $|\dom(\phi)| \leq m$. For $M \subseteq \RM$, define $\dom(M) \defeq \bigcup \set{\dom(\phi) \,:\, \phi \in M}$ and $\mu(M) \defeq \mu(\dom(M))$. For each $x \in V$, there are countably many recoloring moves $\phi$ with $\dom(\phi) \ni x$; it follows from the Luzin--Novikov theorem \cite[Theorem 18.10]{KechrisDST} that $\dom(M)$ is Borel for every Borel set $M \subseteq \RM$.
	
	For a proper $\col{C}$-coloring $f$ of $G$ and $\phi \in \RM$, we define a coloring $f \recolor \phi \colon V \to \col{C}$ by the formula
	\[
	(f \recolor \phi)(x) \,\defeq\, \begin{cases}
	\phi(x) &\text{if } x \in \dom(\phi);\\
	f(x) &\text{otherwise}.
	\end{cases}
	\]
	The coloring $f \recolor \phi$ is the result of \emphd{applying} the recoloring move $\phi$ to $f$, and we say that a recoloring move $\phi$ is \emphd{acceptable} for $f$ if $f \recolor \phi$ is again a proper coloring of $G$. For each color $\alpha \in \col{C}$, let
	\begin{equation}\label{eq:partial}
		\diff_\alpha(f, \phi) \,\defeq\, |\phi^{-1}(\alpha)| - |\dom(\phi) \cap f^{-1}(\alpha)|.
	\end{equation}
	In other words, $\diff_\alpha(f, \phi)$ is the change, between $f$ and $f \recolor \phi$, in the number of vertices of color $\alpha$ among the elements of $\dom(\phi)$. Define the following two \ep{disjoint} sets of colors:
	\begin{equation}\label{eq:inc_dec}
		\col{D}^+(f, \phi) \defeq \set{\alpha \in \col{C} \,:\, \diff_\alpha(f, \phi) > 0} \qquad \text{and} \qquad \col{D}^-(f, \phi) \defeq \set{\alpha \in \col{C} \,:\, \diff_\alpha(f, \phi) < 0}.
	\end{equation}
	Assuming that $f$ is $\mu$-measurable, we say that $\phi$ \emphd{improves $f$ with respect to $\mu$} if
	\begin{enumerate}[label=\ep{\normalfont{}I\arabic*}]
		\item $\phi$ is acceptable for $f$; and
		\item\label{item:I2} there is $\alpha \in \col{D}^+(f, \phi)$ such that for all $\beta \in \col{D}^-(f, \phi)$, we have $\mu(f^{-1}(\alpha)) < \mu(f^{-1}(\beta))$.
	\end{enumerate}
	Informally, \ref{item:I2} states that some ``small'' color class has a net gain of vertices upon the recoloring move $\phi$. Given a set $M \subseteq \RM$, we say that a proper $\col{C}$-coloring $f$ of $G$ is \emphd{$(\mu, M)$-perfect} if
	\[
		\mu(\set{\phi \in M \,:\, \text{$\phi$ improves $f$ with respect to $\mu$}}) \,=\, 0.
	\]
	The main result of this subsection is that $(\mu, \RM_3)$-perfect colorings must be $\mu$-equitable. In other words, if a coloring is {not} $\mu$-equitable, it can be improved by a recoloring move $\phi$ with $|\dom(\phi)| \leq 3$, and, furthermore, such recoloring moves cover a subset of $V$ of positive measure. Our proof of this fact is an adaptation of the proof of the \hyperref[theo:HSz]{Hajnal--Szemer\'edi theorem} from \cite{equit}.
	
	\begin{lemma}\label{lemma:find_move}
		If a $\mu$-measurable proper $\col{C}$-coloring $f$ of $G$ is $(\mu, \RM_3)$-perfect, then it is $\mu$-equitable. 
	\end{lemma}
	\begin{scproof}\stepcounter{ForClaims} \renewcommand{\theForClaims}{\ref{lemma:find_move}}
		Suppose, towards a contradiction, that $f$ is a $\mu$-measurable proper $\col{C}$-coloring of $G$ that is $(\mu, \RM_3)$-perfect but not $\mu$-equitable. After passing to a $\mu$-conull $G$-invariant Borel subset of $V$, we may assume that, in fact, no recoloring move $\phi \in \RM_3$ improves $f$ with respect to $\mu$.
		
		For $\gamma \in \col{C}$, let $V_\gamma \defeq f^{-1}(\gamma)$ be the corresponding color class. Let $a \defeq \min_{\gamma \in \col{C}} \mu(V_\gamma)$ and define \[\col{A} \defeq \set{\alpha \in \col{C} \,:\, \mu(V_\alpha) = a} \qquad \text{and} \qquad \col{B} \defeq \col{C} \setminus \col{A}.\] Since $f$ is not $\mu$-equitable, $\col{B} \neq \0$. Let $b \defeq |\col{B}|^{-1} \sum_{\beta \in \col{B}} \mu(V_\beta)$ and notice that $a < 1/k < b$. Set \[A \defeq f^{-1}(\col{A}) \qquad \text{and} \qquad B \defeq f^{-1}(\col{B}).\]
		
		\begin{claim}\label{claim:move1}
			If $\alpha \in \col{A}$, then every vertex $x \in B$ has a neighbor in $V_\alpha$.
		\end{claim}
		\begin{claimproof}
			Suppose $\beta \in \col{B}$ and $x \in V_\beta$ has no neighbor in $V_\alpha$. Consider the recoloring move $\phi \defeq \set{(x, \alpha)}$ \ep{see Figure~\ref{fig:move1}}. The assumption on $x$ ensures that $\phi$ is acceptable for $f$. Since $\col{D}^+(f, \phi) = \set{\alpha}$, $\col{D}^-(f, \phi) = \set{\beta}$, and $\mu(V_\alpha) < \mu(V_\beta)$,
			we conclude that $\phi$ improves $f$ with respect to $\mu$.
		\end{claimproof}
	
	\begin{figure}[h]
		\centering
		\begin{tikzpicture}
			\draw (0,0) circle [x radius=0.4, y radius=0.8];
			\draw (2,0) circle [x radius=0.4, y radius=0.6];
			
			\node[anchor=south] (x) at (0,0.2) {$x$};
			\filldraw (0,0.2) circle (2pt);
			\draw[-{Stealth[length=1.6mm]},bend right,dashed] (0,0.2) to (1.9,-0.2);
			
			\node at (0,-1.1) {$V_\beta$};
			\node at (2,-0.9) {$V_\alpha$};
		\end{tikzpicture}
		\caption{The recoloring move from Claim~\ref{claim:move1}.}\label{fig:move1}
	\end{figure}
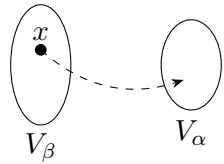
		
		Given $x \in B$ and $y \in A$, we say that $y$ is a \emphd{solo neighbor} of $x$ if $y$ is the unique neighbor of $x$ with the color $f(y)$. \ep{This notion has played a similarly central role in \cite{equit}.} For $y \in A$, let
		\[
			S(y) \,\defeq\, \set{x \in B \,:\, \text{$y$ is a solo neighbor of $x$}}.
		\]
	
		\begin{claim}\label{claim:move2}
			Let $\alpha$, $\alpha' \in \col{A}$ be distinct and let $y \in V_\alpha$. If $S(y) \neq \0$, then $y$ has a neighbor in $V_{\alpha'}$.
		\end{claim}
		\begin{claimproof}
			Take any $x \in S(y)$ and let $\beta \defeq f(x)$. Suppose that $y$ has no neighbor in $V_{\alpha'}$ and consider the recoloring move $\phi \defeq \set{(x, \alpha), (y, \alpha')}$ \ep{see Figure~\ref{fig:move2}}. Since $y$ is not adjacent to any vertex in $V_{\alpha'}$, while the only neighbor of $x$ that is colored $\alpha$ is $y$, $\phi$ is acceptable for $f$. But $\col{D}^+(f, \phi) = \set{\alpha'}$, $\col{D}^-(f, \phi) = \set{\beta}$, and $\mu(V_{\alpha'}) < \mu(V_\beta)$,
			so $\phi$ improves $f$ with respect to $\mu$.
		\end{claimproof}
	
		\begin{figure}[h]
			\centering
			\begin{tikzpicture}
			\draw (0,0) circle [x radius=0.4, y radius=0.8];
			\draw (2,0) circle [x radius=0.4, y radius=0.6];
			\draw (4,0) circle [x radius=0.4, y radius=0.6];
			
			\node[anchor=south] (x) at (0,0.2) {$x$};
			\filldraw (0,0.2) circle (2pt);
			\draw[-{Stealth[length=1.6mm]},bend right,dashed] (0,0.2) to (1.9,-0.2);
			
			\node[anchor=south] (y) at (2,0.0) {$y$};
			\filldraw (2,0.0) circle (2pt);
			\draw[-{Stealth[length=1.6mm]},bend left,dashed] (2,0.0) to (4,0.3);
			
			\draw[thick] (0,0.2) to (2,0.0); 
			
			\node at (0,-1.1) {$V_\beta$};
			\node at (2,-0.9) {$V_\alpha$};
			\node at (4,-0.9) {$V_{\alpha'}$};
			\end{tikzpicture}
			\caption{The recoloring move from Claim~\ref{claim:move2}.}\label{fig:move2}
		\end{figure}
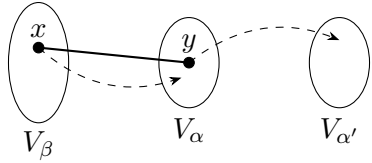
	
		\begin{claim}\label{claim:move3}
			If $y \in A$, then the induced subgraph $G[S(y)]$ is a clique.
		\end{claim}
		\begin{claimproof}
			Suppose, towards a contradiction, that $G[S(y)]$ is not a clique. This means that there are two distinct non-adjacent vertices $x$, $x' \in S(y)$. Let $\alpha \defeq f(y)$, $\beta \defeq f(x)$, and $\beta' \defeq f(x')$. Observe that there is a color $\gamma \neq \alpha$ such that \begin{equation}\label{eq:extra_color}
				(N_G(y) \cap V_\gamma) \setminus \set{x, x'} \,=\, \0.
			\end{equation}
			Otherwise, $y$ would have, in addition to $x$ and $x'$, at least one neighbor in every color class $V_\gamma$ except for $V_\alpha$, which would yield $\deg_G(y) \geq 2 + (k-1) = k+1 > \Delta$. Now, fix any $\gamma \neq \alpha$ satisfying \eqref{eq:extra_color} and consider the recoloring move $\phi \defeq \set{(x, \alpha),(x', \alpha),(y, \gamma)}$ \ep{see Figure~\ref{fig:move3}}. The choice of $x$, $x'$, and $\gamma$ implies that $\phi$ is acceptable for $f$. Since $\alpha \in \col{D}^+(f, \phi)$, $\col{D}^-(f, \phi) \subseteq \set{\beta, \beta'}$, and $\mu(V_\alpha) < \mu(V_\beta)$, $\mu(V_{\beta'})$,
			we conclude that $\phi$ improves $f$ with respect to $\mu$.
		\end{claimproof}
	
		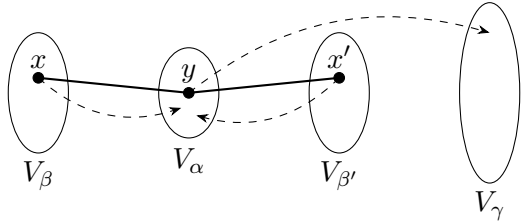
\begin{figure}[h]
			\centering
			\begin{tikzpicture}
			\draw (0,0) circle [x radius=0.4, y radius=0.8];
			\draw (2,0) circle [x radius=0.4, y radius=0.6];
			\draw (4,0) circle [x radius=0.4, y radius=0.8];
			\draw (6,0) circle [x radius=0.4, y radius=1.2];
			
			\node[anchor=south] (x) at (0,0.2) {$x$};
			\filldraw (0,0.2) circle (2pt);
			\draw[-{Stealth[length=1.6mm]},bend right,dashed] (0,0.2) to (1.9,-0.2);
			
			\node[anchor=south] (y) at (2,0.0) {$y$};
			\filldraw (2,0.0) circle (2pt);
			\draw[-{Stealth[length=1.6mm]},dashed] (2,0.0) to[bend left=27] (6,0.8);
			
			\node[anchor=south] at (4,0.2) {$x'$};
			\filldraw (4,0.2) circle (2pt);
			\draw[-{Stealth[length=1.6mm]},bend left,dashed] (4,0.2) to (2.1,-0.3);
			
			\draw[thick] (0,0.2) to (2,0.0); 
			\draw[thick] (2,0.0) to (4,0.2);
			
			\node at (0,-1.1) {$V_\beta$};
			\node at (2,-0.9) {$V_\alpha$};
			\node at (4,-1.1) {$V_{\beta'}$};
			\node at (6,-1.5) {$V_\gamma$};
			\end{tikzpicture}
			\caption{The recoloring move from Claim~\ref{claim:move3}.}\label{fig:move3}
		\end{figure}
	
		\begin{claim}\label{claim:B_small}
			We have $\Delta \geq 2|\col{B}|+1$.
		\end{claim}
		\begin{claimproof}
			Fix an arbitrary color $\alpha \in \col{A}$. Let $X$ be the set of all vertices $x \in B$ that have a solo neighbor in $V_\alpha$ and let $Y$ be the set of all vertices $y \in V_\alpha$ that have a neighbor in $X$. By Claim~\ref{claim:move2}, every vertex $y \in Y$ has at least $|\col{A}| - 1 = k - |\col{B}| - 1$ neighbors in $A$, and thus
			\begin{equation}\label{eq:deg_bound}
				|N_G(y) \cap B| \,\leq\, \Delta - k + |\col{B}| + 1.
			\end{equation}
			This immediately yields
			\begin{equation}\label{eq:upper}
			\int_{V_{\alpha}} |N_G(y) \cap B| \,\Diff \mu(y) \,\leq\, \Delta \mu(V_\alpha) \,-\, (k - |\col{B}| - 1) \mu(Y).
			\end{equation}
			On the other hand, we have
			\begin{align*}
				\mu(X) \,&=\, \int_{X} 1 \,\Diff \mu(x) \\
				[\text{by the definition of $X$ and $Y$}]\qquad&=\, \int_{X} |N_G(x) \cap Y| \,\Diff\mu(x) \\
				[\text{$\mu$ is $G$-invariant}]\qquad&=\, \int_{Y} |N_G(y) \cap X| \,\Diff\mu(y) \\
				[\text{by \eqref{eq:deg_bound}}]\qquad&\leq\, (\Delta - k + |\col{B}| + 1) \mu(Y).
			\end{align*}
			This, together with Claim~\ref{claim:move1}, gives
			\begin{equation}\label{eq:lower}
				\int_{B} |N_G(x) \cap V_\alpha| \,\Diff \mu(x) \,\geq\, 2\mu(B) \,-\, \mu(X) \,\geq\, 2\mu(B) \,-\, (\Delta - k + |\col{B}| + 1) \mu(Y).
			\end{equation}
			Since $\mu$ is $G$-invariant, we can combine \eqref{eq:lower} and \eqref{eq:upper} to get
			\[
				\Delta \mu(V_\alpha) \,-\, (k - |\col{B}| - 1) \mu(Y) \,\geq\, 2\mu(B) \,-\, (\Delta - k + |\col{B}| + 1) \mu(Y).
			\]
			Recalling that $\mu(V_\alpha) = a$ and $\mu(B) = |\col{B}|b$, we can rewrite the last inequality as
			\[
				\Delta a \,-\, (k-|\col{B}|-1)\mu(Y) \,\geq\, 2|\col{B}|b \,-\, (\Delta-k+|\col{B}|+1)\mu(Y).
			\]
			After moving all the terms to one side and using that $b > a$, we obtain
			\begin{equation*}
				\Delta \,-\, 2|\col{B}| \,+\, (\Delta -2k + 2 |\col{B}| +2)\frac{\mu(Y)}{b} \,>\, 0.
			\end{equation*}
			Since $0 \leq \mu(Y)/b < \mu(Y)/a \leq 1$, at least one of the following two inequalities holds:
			\[
				\Delta - 2|\col{B}| \,>\, 0 \qquad \text{or} \qquad \Delta \,-\, 2|\col{B}| \,+\, (\Delta -2k + 2 |\col{B}| +2) \,=\, 2\Delta -2k +2 \,>\, 0.
			\]
			The second of these inequalities necessarily fails, since $2\Delta - 2k + 2 \leq 2\Delta - 2(\Delta + 1) + 2 = 0$. Hence, we must have $\Delta - 2|\col{B}| > 0$, and thus $\Delta \geq 2|\col{B}|+1$, as desired.
		\end{claimproof}
	
		We are now ready for the final stage of the argument. At least one of the color classes $V_\beta$, $\beta \in \col{B}$, has measure at least $b$, so by taking $J$ to be any such color class and applying Corollary~\ref{corl:max_indep} to the induced subgraph $G[B]$, we obtain a Borel maximal $G$-independent subset $I \subseteq B$ such that $\mu(I) \geq b$. For each $x \in I$, let $\Sigma(x)$ be the set of all solo neighbors of $x$. Then, on the one hand,
		\[
			|N_G(x) \cap A| \,\geq\, 2|\col{A}| - |\Sigma(x)| \,=\, 2k - 2|\col{B}| - |\Sigma(x)|,
		\]
		while, on the other hand, $|N_G(x) \cap A| \leq \Delta - |N_G(x) \cap B|$, which yields
		\begin{equation}\label{eq:Sigma}
			|\Sigma(x)| \,\geq\, 2k - 2|\col{B}| - \Delta + |N_G(x) \cap B|.
		\end{equation}
		Since $I$ is $G$-independent, Claim~\ref{claim:move3} implies that the sets $\Sigma(x)$, $x \in I$, are pairwise disjoint. Using the fact that $\bigcup_{x \in I} \Sigma(x) \subseteq A$ and the $G$-invariance of $\mu$, we conclude that
		\[
			\mu(A) \,\geq\, \int_{I} |\Sigma(x)| \,\Diff\mu(x).
		\]
		From \eqref{eq:Sigma}, it follows that
		\[
			\int_I |\Sigma(x)| \,\Diff\mu(x) \,\geq\, (2k - 2|\col{B}| - \Delta) \mu(I) \,+\, \int_I |N_G(x) \cap B| \,\Diff \mu(x). 
		\]
		Since $I$ is a {maximal} $G$-independent subset of $B$, every vertex $z \in B \setminus I$ has a neighbor in $I$, and thus
		\[
			\int_I |N_G(x) \cap B| \,\Diff \mu(x) \,=\, \int_{B \setminus I} |N_G(z) \cap I| \,\Diff \mu(z) \,\geq\, \mu(B \setminus I) \,=\, \mu(B) - \mu(I).
		\]
		To summarize, we have
		\begin{equation}\label{eq:A_large}
			\mu(A) \,\geq\, (2k-2|\col{B}|-\Delta-1) \mu(I) \,+\, \mu(B).
		\end{equation}
		Claim~\ref{claim:B_small} and the inequality $k \geq \Delta+1$ yield $2k-2|\col{B}|-\Delta - 1 > 0$. Hence, since $\mu(I) \geq b$, the right-hand side of \eqref{eq:A_large} could only decrease if we replace $\mu(I)$ by $b$, so
		\[
			\mu(A) \,\geq\, (2k-2|\col{B}|-\Delta - 1) b + \mu(B).
		\]
		Using that $\mu(B) = |\col{B}|b$, $\mu(A) = |\col{A}|a = (k - |\col{B}|)a$, and $b > a$, we obtain
		\[
			k - |\col{B}| \,>\, (2k-2|\col{B}|-\Delta - 1) + |\col{B}|,
		\]
		which after simplifying becomes $k < \Delta + 1$. This is a contradiction, and Lemma~\ref{lemma:find_move} is proved.
	\end{scproof}

	Roughly speaking, our strategy now is to start with an arbitrary proper coloring $f_0$ and repeatedly apply to $f_0$ recoloring moves that improve it, producing an infinite sequence of colorings $(f_n)_{n=0}^\infty$. We then intend to show that this sequence converges to some ``limit'' coloring $f_\infty$ that cannot be improved by a recoloring move anymore \ep{at least on a $\mu$-conull set}, and hence, by Lemma~\ref{lemma:find_move}, $f_\infty$ must be $\mu$-equitable. For this approach to work, we must argue that the ``limit'' coloring $f_\infty$ actually exists. This will be done using the Borel--Cantelli lemma: we will prove that $\sum_{n = 0}^\infty \dist_\mu(f_n, f_{n+1}) < \infty$, which implies that the pointwise limit $\lim_{n \to \infty} f_n(x)$ exists for $\mu$-almost every $x \in V$. To obtain an upper bound on $\sum_{n = 0}^\infty \dist_\mu(f_n, f_{n+1})$, we will require a certain technical result \ep{namely Lemma~\ref{lemma:limit}} that is proven in the next subsection. At this point we should note that our argument would be quite a bit simpler \ep{and we would have no need of Lemma~\ref{lemma:limit}} if it were possible to control the recoloring process by some numerical parameter, for instance, if we could ensure that $\disc_\mu(f_{n+1}) < \disc_\mu(f_n)$. Unfortunately, this is not the case; the culprit is the recoloring move from Claim~\ref{claim:move3} \ep{see Figure~\ref{fig:move3}}, during which a vertex may be added to a color class \ep{namely $V_\gamma$} that is already ``too large.''

	\subsection{Comparing distributions}\label{subsec:distrib}
	
	\begin{assum*}[for \S\ref{subsec:distrib}]
		Fix a finite set of colors $\col{C}$ of size $k \geq 1$.
	\end{assum*}

	As usual, for a function $\omega \colon \col{C} \to \R$, we write $\|\omega\|_1 \defeq \sum_{\alpha \in \col{C}} |\omega(\alpha)|$. A \emphd{distribution \ep{on $\col{C}$}} is a function $\omega \colon \col{C} \to [0,1]$ such that $\|\omega\|_1 = 1$. The \emphd{discrepancy} of a distribution $\omega$ is \[\disc(\omega) \,\defeq\, \max_{\alpha \in \col{C}} \left|\omega(\alpha) - 1/k\right|.\] Given a pair of distributions $\omega$ and $\eta$, we define the following two \ep{disjoint} sets of colors:
	\[
		\col{D}^+(\omega, \eta) \defeq \set{\alpha \in \col{C} \,:\, \eta(\alpha) > \omega(\alpha)} \qquad \text{and} \qquad \col{D}^-(\omega, \eta) \defeq \set{\alpha \in \col{C} \,:\, \eta(\alpha) < \omega(\alpha)}.
	\]
	Observe the similarity between this definition and \eqref{eq:inc_dec}. We say that $\eta$ is \emphd{more equitable} than $\omega$, in symbols $\omega \vartriangleleft \eta$, if there is a color $\alpha \in \col{D}^+(\omega, \eta)$ such that for all $\beta \in \col{D}^-(\omega, \eta)$, we have $\eta(\alpha) \leq \eta(\beta)$, and we write $\omega \trianglelefteq \eta$ to mean that $\omega \vartriangleleft \eta$ or $\omega = \eta$.
	Notice that the relation $\vartriangleleft$ is antisymmetric and irreflexive; however, it is not transitive when $k \geq 3$. Nevertheless, the {transitive closure} of $\vartriangleleft$ is a strict partial order on distributions, which justifies our use of the order-like symbol ``$\vartriangleleft$.''\footnote{Since we will not use this fact explicitly, its proof is not included; but it follows readily from Claim~\ref{claim:initial_sums}.} We also point out that if $\eta(\alpha) = 1/k$ for all $\alpha \in \col{C}$ \ep{i.e., if $\eta$ is the \emphd{uniform distribution}}, then $\eta$ is more equitable than all other distributions $\omega$. 
	
	
	\begin{lemma}\label{lemma:limit}
		Let $(\omega_n)_{n=0}^\infty$ be a sequence of distributions on $\col{C}$ such that for all $n \in \N$, $\omega_n \trianglelefteq \omega_{n+1}$. Suppose that there is a real number $A \geq 1$ such that for all $n \in \N$ and $\alpha \in \col{D}^+(\omega_n, \omega_{n+1})$,
		\[
			\|\omega_{n+1} - \omega_{n}\|_1 \,\leq\, A \cdot \left(\omega_{n+1}(\alpha) - \omega_n(\alpha)\right).
		\]
		Then
		\[
			\sum_{n=0}^\infty \|\omega_{n+1} - \omega_{n}\|_1 \,\leq\, \frac{(1+A)^{k+1}}{A} \cdot \disc(\omega_0) \,<\, \infty. 
		\]
	\end{lemma}
	\begin{scproof}\stepcounter{ForClaims} \renewcommand{\theForClaims}{\ref{lemma:limit}}
		For each distribution $\omega$ on $\col{C}$, let $\omega^\ast \colon \set{1, \ldots, k} \to [0,1]$ denote the mapping obtained by putting the values of $\omega$ in non-decreasing order; that is, the lists
		\[
			\omega(\alpha), \, \alpha \in \col{C} \qquad \text{and} \qquad \omega^\ast(1), \, \omega^\ast(2), \, \ldots, \, \omega^\ast(k)
		\]
		contain the same elements, and also $\omega^\ast(1) \leq \omega^\ast(2) \leq \ldots \leq \omega^\ast(k)$.
		
		\begin{claim}\label{claim:ast}
			For all distributions $\omega$ and $\eta$, we have $\|\omega^\ast - \eta^\ast\|_1 \leq \|\omega - \eta\|_1$.
		\end{claim}
		\begin{claimproof}
			For each $\alpha \in \col{C}$, let $I_\alpha$ be the open interval with endpoints $\omega(\alpha)$ and $\eta(\alpha)$. Similarly, for each $1 \leq i \leq k$, let $J_i$ be the open interval with endpoints $\omega^\ast(i)$ and $\eta^\ast(i)$. Then
			\begin{equation}\label{eq:lengths}
				\|\omega - \eta\|_1 \,=\, \sum_{\alpha \in \col{C}} \mathrm{length}(I_\alpha) \qquad \text{and} \qquad \|\omega^\ast - \eta^\ast\|_1 \,=\, \sum_{i=1}^k \mathrm{length}(J_i).
			\end{equation}
			Let $S \defeq \set{\omega(\alpha),\, \eta(\alpha) \,:\, \alpha \in \col{C}}$. Note that the set $S$ is finite. We shall argue that for each $x \in \R\setminus S$, the number of colors $\alpha$ with $x \in I_\alpha$ is not less than the number of indices $i$ with $x \in J_i$; in view of \eqref{eq:lengths}, this immediately yields the claim. Take any $x \in \R\setminus S$ and define
			\[
				\tau(x) \,\defeq\, |\set{\alpha \in \col{C} \,:\, \omega(\alpha) < x}| \quad \text{and} \quad \sigma(x) \,\defeq\, |\set{\alpha \in \col{C} \,:\, \eta(\alpha) < x}|.
			\]
			Then $x$ is in {precisely} $|\tau(x) - \sigma(x)|$ of the intervals $J_i$, $1 \leq i \leq k$. On the other hand, if $\tau(x) \geq \sigma(x)$, then there are at least $\tau(x) - \sigma(x)$ colors $\alpha$ with $\omega(\alpha) < x < \eta(\alpha)$, while if $\sigma(x) \geq \tau(x)$, then there are at least $\sigma(x) - \tau(x)$ colors $\alpha$ with $\eta(\alpha) < x < \omega(\alpha)$. In either case, $x$ belongs to {at least} $|\tau(x) - \sigma(x)|$ of the intervals $I_\alpha$, $\alpha \in \col{C}$, and hence we are done.
		\end{claimproof}
		
		The proof of Lemma~\ref{lemma:limit} rests on the following key observation:
		
		\begin{claim}\label{claim:initial_sums}
			Suppose that $\omega$ and $\eta$ are distributions on $\col{C}$ and $\omega \vartriangleleft \eta$. Then there is an index $\ell$ such that $\omega^\ast(i) \leq \eta^\ast(i)$ for all $1 \leq i \leq \ell$ and there is a color $\alpha \in \col{D}^+(\omega, \eta)$ with
			\begin{equation}\label{eq:sum_bound}
				\sum_{i=1}^\ell \eta^\ast(i) \,-\, \sum_{i=1}^\ell \omega^\ast(i) \,\geq\, \eta(\alpha) - \omega(\alpha).
			\end{equation}
		\end{claim}
		\begin{claimproof}
			Let $\alpha \in \col{D}^+(\omega, \eta)$ be such that for all $\beta \in \col{D}^-(\omega, \eta)$, we have $\eta(\alpha) \leq \eta(\beta)$, and let $\ell$ be the least index such that $\eta^\ast(\ell) = \eta(\alpha)$. We claim that this choice of $\ell$ and $\alpha$ works.
			
			To begin with, fix a bijection $\set{1, \ldots, k} \to \col{C} \colon i \mapsto \alpha_i$ such that $\eta^\ast(i) = \eta(\alpha_i)$ for all $1 \leq i \leq k$, and $\alpha_\ell = \alpha$. Consider any $1 \leq i \leq \ell$. We have $\eta(\alpha_i) = \eta^\ast(i) \leq \eta^\ast(\ell)$, where equality holds if and only if $i = \ell$. By the choice of $\alpha$, this yields $\alpha_i \not \in \col{D}^-(\omega, \eta)$, and therefore $\omega(\alpha_i) \leq \eta(\alpha_i) = \eta^\ast(i)$. Thus, there are at least $i$ colors---namely $\alpha_1$, \ldots, $\alpha_i$---at which the value of $\omega$ does not exceed $\eta^\ast(i)$. But this precisely means that $\eta^\ast(i) \geq \omega^\ast(i)$, as desired.
			
			Next, we prove \eqref{eq:sum_bound}. Let $1 \leq t \leq \ell$ be the largest index such that $\omega^\ast(t) \leq \omega(\alpha)$ \ep{such $t$ exists as $\omega^\ast(1) \leq \omega(\alpha)$}. We claim that if $t < s \leq \ell$, then $\omega^\ast(s) \leq \eta^\ast(s-1)$. Indeed, suppose, towards a contradiction, that $\omega^\ast(s) > \eta^\ast(s-1)$. Then $\omega^\ast(s-1) \leq \eta^\ast(s-1) < \omega^\ast(s)$, meaning that the only colors at which the value of $\omega$ is at most $\eta^\ast(s-1)$ are $\alpha_1$, \ldots, $\alpha_{s-1}$. Since $\alpha = \alpha_\ell$ is not among them, $\omega(\alpha) \geq \omega^\ast(s)$, contradicting the choice of $t$ and the fact that $s > t$. Now we can write
			\begin{align*}
				\sum_{i=1}^\ell \eta^\ast(i) \,-\, \sum_{i=1}^\ell \omega^\ast(i) \,&=\, \eta^\ast(\ell) \,+\, \sum_{s=t+1}^\ell \left(\eta^\ast(s-1) - \omega^\ast(s)\right) \\
				&\qquad\qquad-\, \omega^\ast(t) \,+\, \sum_{i = 1}^{t-1} \left(\eta^\ast(i) - \omega^\ast(i)\right) \\
				[\text{the other summands are nonnegative}]\qquad &\geq\, \eta^\ast(\ell) - \omega^\ast(t) \\
				[\text{since $\eta^\ast(\ell) = \eta(\alpha)$ and $\omega^\ast(t) \leq \omega(\alpha)$}]\qquad &\geq\, \eta(\alpha) - \omega(\alpha),
			\end{align*}
			as desired.
		\end{claimproof}
	
		Now let $(\omega_n)_{n=0}^\infty$ be as in the statement of Lemma~\ref{lemma:limit}. Define $S_n(\ell) \defeq \sum_{i=1}^\ell \omega^\ast_n(i)$ for each $n \in \N$ and $1 \leq \ell \leq k$. Let $R \defeq \set{n \in \N\,:\, \omega_n \neq \omega_{n+1}}$. Claim~\ref{claim:initial_sums} shows that for every $n \in R$, there exist an index $\ell_n$ and a color $\alpha_n \in \col{D}^+(\omega, \eta)$ such that $\omega^\ast_n(i) \leq \omega^\ast_{n+1}(i)$ for all $1 \leq i \leq \ell_n$, and also
		\begin{equation}\label{eq:bound}
			S_{n+1}(\ell_n) - S_n(\ell_n) \,\geq\, \omega_{n+1}(\alpha_n) - \omega_n(\alpha_n) \,\geq\, A^{-1} \cdot \|\omega_{n+1} - \omega_n\|_1.
		\end{equation}
		For each $\ell$, let $R_\ell \defeq \set{n \in R \,:\, \ell_n = \ell}$ and $R_{<\ell} \defeq \set{n \in R \,:\, \ell_n < \ell} = \bigcup_{i=1}^{\ell-1} R_i$. We will show that
		\begin{equation}\label{eq:sum_of_sums}
			\sum_{n \in R_\ell} \left(S_{n+1}(\ell) - S_n(\ell)\right) \,\leq\, \frac{(1+A)^\ell - 1}{A} \cdot \disc(\omega_0).
		\end{equation}
		Once inequality \eqref{eq:sum_of_sums} is proved, we obtain, using \eqref{eq:bound}, that 
		\begin{align*}
			\sum_{n=0}^\infty \|\omega_{n+1} - \omega_n\|_1 \,&\leq\, A \cdot \sum_{n\in R} \left(S_{n+1}(\ell_n) - S_n(\ell_n)\right) \\
			&=\, A \cdot \sum_{\ell=1}^k \sum_{n \in R_\ell} \left(S_{n+1}(\ell) - S_n(\ell)\right) \\
			[\text{by \eqref{eq:sum_of_sums}}]\qquad&\leq\, \sum_{\ell=1}^k \left((1+A)^\ell - 1\right) \cdot \disc(\omega_0) \,=\, \frac{(1+A)^{k+1} - (k+1)A - 1}{A}\cdot \disc(\omega_0) \\
			&\leq\, \frac{(1+A)^{k+1}}{A} \cdot \disc(\omega_0),
		\end{align*}
		as desired. To prove \eqref{eq:sum_of_sums}, we use induction on $\ell$, so suppose that \eqref{eq:sum_of_sums} holds with $\ell$ replaced by any $1 \leq i < \ell$. Notice that if $\ell_n \geq \ell$, then $S_{n+1}(\ell) - S_n(\ell) \geq 0$. Hence, for each $m \in \N$, 
		\begin{align*}
		    \ell \cdot \disc(\omega_0)  \,\geq\, \frac{\ell}{k} \,-\, S_0(\ell) \,&\geq\, S_{m+1}(\ell) \,-\, S_0(\ell) \,=\, \sum_{n = 0}^{m} \left(S_{n+1}(\ell) - S_n(\ell)\right) \\
		    &\geq\, \sum_{\substack{n \in R_\ell, \\ n \leq m}} \left(S_{n+1}(\ell) - S_n(\ell)\right) \,-\, \sum_{\substack{n \in R_{< \ell}, \\ n \leq m}} \left|S_{n+1}(\ell) - S_n(\ell)\right|.
		\end{align*}
		As $m$ here is arbitrary, we conclude that
		\begin{align*}
			\sum_{n \in R_\ell} \left(S_{n+1}(\ell) - S_n(\ell)\right) 
			\,\leq\, \ell \cdot \disc(\omega_0) \,+\, \sum_{n \in R_{< \ell}} \left|S_{n+1}(\ell) - S_n(\ell)\right|.
		\end{align*}
		We can now write the following chain of inequalities: 
		\begin{align*}
			\sum_{n \in R_{< \ell}} \left|S_{n+1}(\ell) - S_n(\ell)\right| \,&\leq\, \sum_{n \in R_{< \ell}} \left\|\omega^\ast_{n+1} - \omega^\ast_n \right\|_1 \\
			[\text{by Claim~\ref{claim:ast}}]\qquad&\leq\, \sum_{n \in R_{< \ell}} \left\|\omega_{n+1} - \omega_n \right\|_1 \\
			[\text{by \eqref{eq:bound}}]\qquad&\leq\, A \cdot \sum_{i=1}^{\ell-1} \sum_{n \in R_i} \left(S_{n+1}(i) - S_n(i)\right) \\
			[\text{by the inductive hypothesis}]\qquad&\leq\, \sum_{i = 1}^{\ell-1} \left((1+A)^i - 1\right) \cdot \disc(\omega_0).
		\end{align*}
		Therefore,
		\[
			\sum_{n \in R_\ell} \left(S_{n+1}(\ell) - S_n(\ell)\right) \,\leq\, \ell \cdot \disc(\omega_0)+ \sum_{i = 1}^{\ell-1} \left((1+A)^i - 1\right) \cdot \disc(\omega_0) \,=\, \frac{(1+A)^{\ell}-1}{A} \cdot \disc(\omega_0),
		\]
		and the proof is complete.
	\end{scproof}

	\subsection{Perfecting a coloring}\label{subsec:elimination}
	
	\begin{assum*}[for \S\ref{subsec:elimination}]
		Fix an aperiodic Borel graph $G$ of finite maximum degree $\Delta$ with vertex set $V$ and edge set $E$ and a finite set of colors $\col{C}$ of size $k \geq \Delta + 1$.
	\end{assum*}

	Let $f \colon V \to \col{C}$ be a Borel coloring. Given a probability measure $\mu$ on $V$, $f$ gives rise to the push\-forward distribution $f_\ast(\mu)$ on $\col{C}$ defined by $f_\ast(\mu)(\alpha) \defeq \mu(f^{-1}(\alpha))$. Note that $\disc(f_\ast(\mu)) = \disc_\mu(f)$.
	
	\begin{lemma}\label{lemma:dist1}
		If $\mu$ is a probability measure on $V$ and $f$, $g \colon V \to \col{C}$ are $\mu$-measurable, then
		\[
			\|f_\ast(\mu)-g_\ast(\mu)\|_1 \,\leq\, 2 \cdot \dist_\mu(f,g).
		\]
	\end{lemma}
	\begin{scproof}
		For $\alpha$, $\beta \in \col{C}$, let $V_{\alpha\beta} \defeq f^{-1}(\alpha) \cap g^{-1}(\beta)$. Then, for each fixed $\alpha \in \col{C}$,
		\[
			|\mu(f^{-1}(\alpha)) - \mu(g^{-1}(\alpha))| \,=\, \left|\sum_{\beta \in \col{C}} \mu(V_{\alpha \beta}) \,-\, \sum_{\beta \in \col{C}} \mu(V_{\beta \alpha}) \right| \,\leq\, \sum_{\substack{\beta \in \col{C}:\\\beta \neq \alpha}} \left(\mu(V_{\alpha \beta}) + \mu(V_{\beta \alpha})\right).
		\]
		and therefore,
		\[
			\|f_\ast(\mu)-g_\ast(\mu)\|_1 \,\leq\, \sum_{\alpha \in \col{C}} \sum_{\substack{\beta \in \col{C}:\\\beta \neq \alpha}} \left(\mu(V_{\alpha \beta}) + \mu(V_{\beta \alpha})\right) \,=\, 2 \cdot \sum_{\substack{\alpha, \beta \in \col{C}:\\\alpha \neq \beta}} \mu(V_{\alpha \beta}) \,=\, 2 \cdot \dist_\mu(f,g).\qedhere
		\]
	\end{scproof}
	
	Starting with a Borel proper $\col{C}$-coloring $f$, we wish to apply recoloring moves in order to make the distribution $f_\ast(\mu)$ more equitable. Since applying a single recoloring move changes $f$ only on a finite set of vertices, we have to be able to apply infinitely many recoloring moves simultaneously, but we must take care that the moves do not interfere with each other. Say that a set $M \subseteq \RM$ is \emphd{$G$-separated} if for every two distinct recoloring moves $\phi$, $\phi' \in M$,
	\begin{enumerate}[label=\ep{\normalfont{}S\arabic*}]
		\item\label{item:disjoint} $\dom(\phi) \cap \dom(\phi') = \0$; and
		\item\label{item:no_edges} there are no edges in $G$ between $\dom(\phi)$ and $\dom(\phi')$.
	\end{enumerate}
	For a Borel $G$-separated set $M \subseteq \RM$, define a \ep{Borel} coloring $f \recolor M \colon V \to \col{C}$ via
	\[
	(f \recolor M)(x) \,\defeq\, \begin{cases}
	\phi(x) &\text{if } \phi \in M \text{ and } x \in \dom(\phi);\\
	f(x) &\text{if } x \not \in \dom(M).
	\end{cases}
	\]
	In other words, $f \recolor M$ is the result of simultaneously applying every recoloring move $\phi \in M$. This is well-defined by \ref{item:disjoint}. Notice that if every recoloring move $\phi \in M$ is acceptable for $f$, then $f \recolor M$ is a proper coloring of $G$ by \ref{item:no_edges}.
	
	To better control the properties of $f \recolor M$, it is convenient to assume that the recoloring moves in $M$ are somewhat similar to each other. Specifically, given a Borel proper $\col{C}$-coloring $f$ of $G$, an integer $m \in \N^+$, and disjoint nonempty sets $\col{D}^+$, $\col{D}^- \subseteq \col{C}$, let $\RM_m(f; \col{D}^+, \col{D}^-)$ denote the set of all recoloring moves $\phi \in \RM_m$ such that $\phi$ is acceptable for $f$ and
	\[
		\col{D}^+(f, \phi) = \col{D}^+ \qquad \text{and} \qquad \col{D}^-(f, \phi) = \col{D}^-.
	\]
	Note that the set $\RM_m(f; \col{D}^+, \col{D}^-)$ is Borel. If $M \subseteq \RM_m(f;\col{D}^+, \col{D}^-)$ is a Borel $G$-separated subset and $\mu$ is a $G$-invariant probability measure on $V$, then either $\mu(M) = 0$, or else,
	\[
		\col{D}^+(f_\ast(\mu), (f \recolor M)_\ast(\mu)) = \col{D}^+ \qquad \text{and} \qquad \col{D}^-(f_\ast(\mu), (f \recolor M)_\ast(\mu)) = \col{D}^-.
	\]
	
	\begin{lemma}\label{lemma:dist2}
		Let $f$ be a Borel proper $\col{C}$-coloring of $G$. Fix $m \in \N^+$ and disjoint nonempty sets $\col{D}^+$, $\col{D}^- \subseteq \col{C}$. If $M \subseteq \RM_m(f; \col{D}^+, \col{D}^-)$ is a Borel $G$-separated set, then for every $G$-invariant probability measure $\mu$:
		\begin{enumerate}[label=\ep{\normalfont\arabic*}]
			\item $\dist_\mu(f,f \recolor M) \leq m \cdot \|f_\ast(\mu)-(f \recolor M)_\ast(\mu)\|_1$;
			\item for all $\alpha \in \col{D}^+$, we have $\|f_\ast(\mu)-(f \recolor M)_\ast(\mu)\|_1 \leq 2m \cdot (\mu((f \recolor M)^{-1}(\alpha)) - \mu(f^{-1}(\alpha)))$.
		\end{enumerate}
	\end{lemma}
	\begin{scproof}
		Since each $\phi \in M$ has finite domain, there is a Borel function $p \colon M \to V$ such that $p(\phi) \in \dom(\phi)$ for every $\phi \in M$ \ep{see \S\ref{subsec:finite}}. Set $P \defeq \im(p)$. Note that $P$ is Borel, as it is the image of a Borel set under an injective Borel function \cite[Corollary 15.2]{KechrisDST}. Since for each $\phi \in M$, $|\dom(\phi) \cap P| = 1$ and $|\dom(\phi)| \leq m$, and since $\mu$ is $G$-invariant, we conclude that $\mu(P) \geq \mu(M)/m$. For each $x \in P$, let $\phi_x \in M$ be the unique recoloring move such that $p(\phi_x) = x$. The $G$-invariance of $\mu$ yields that for each $\alpha \in \col{C}$,
		\[
			\mu\left((f \recolor M)^{-1}(\alpha)\right) \,-\, \mu(f^{-1}(\alpha)) \,=\, \int_P \diff_\alpha(f, \phi_x) \,\Diff\mu(x),
		\]
		where $\diff_\alpha(f, \phi)$ is defined in \eqref{eq:partial}. Crucially, $M \subseteq \RM_m(f; \col{D}^+, \col{D}^-)$, which means that $\diff_\alpha(f, \phi) \geq 1$ for all $\alpha \in \col{D}^+$ and $\phi \in M$. Hence, for any $\alpha \in \col{D}^+$, we can write
		\begin{align*}
			\|f_\ast(\mu)-(f \recolor M)_\ast(\mu)\|_1 \,&\geq\, \mu\left((f \recolor M)^{-1}(\alpha)\right) \,-\, \mu(f^{-1}(\alpha)) \\
			&\geq\, \int_P 1 \,\Diff\mu(x)
			\,=\, \mu(P) \,\geq\, \frac{\mu(M)}{m} \\
			&\geq\, \frac{\dist_\mu(f, f \recolor M)}{m}\\
			[\text{by Lemma~\ref{lemma:dist1}}]\qquad&\geq\, \frac{1}{2m}\cdot \|f_\ast(\mu)-(f \recolor M)_\ast(\mu)\|_1.\qedhere
		\end{align*}
	\end{scproof}
	
	The next lemma is the main result of this subsection.
	
	\begin{lemma}\label{lemma:step}
		Let $f$ be a Borel proper $\col{C}$-coloring of $G$. Fix an integer $m \in \N^+$ and disjoint nonempty sets $\col{D}^+$, $\col{D}^- \subseteq \col{C}$. Let $M \subseteq \RM_m(f; \col{D}^+, \col{D}^-)$ be a Borel $G$-separated set.
		
		\begin{enumerate}[label=\ep{\itshape\alph*}]
		    \item\label{item:step_a} For every 
		$G$-invariant probability measure $\mu$, there is a Borel proper $\col{C}$-coloring $g$ such that: 
		\begin{enumerate}[label=\ep{\normalfont{}P\arabic*}]
			\item\label{item:dist_norm}\label{item:first} $\dist_\mu(f,g) \leq m \cdot \|f_\ast(\mu)-g_\ast(\mu)\|_1$;
			\item\label{item:norm_step} for all $\alpha \in \col{D}^+(f_\ast(\mu), g_\ast(\mu))$, we have $\|f_\ast(\mu)-g_\ast(\mu)\|_1 \leq 2m \cdot (\mu(g^{-1}(\alpha)) - \mu(f^{-1}(\alpha)))$;
			\item $f_\ast(\mu) \trianglelefteq g_\ast(\mu)$; and
			\item\label{item:last} $g$ is $(\mu, M)$-perfect.
		\end{enumerate}
		\item\label{item:step_b} Furthermore, $G$ admits a Borel proper $\col{C}$-coloring $g$ such that statements \ref{item:first}--\ref{item:last} hold for every ergodic $G$-invariant probability measure $\mu$ simultaneously.
		\end{enumerate}
	\end{lemma}
	\begin{scproof}\stepcounter{ForClaims} \renewcommand{\theForClaims}{\ref{lemma:step}}
		Since $G$ is aperiodic, every $G$-invariant probability measure is atomless, so if $M$ is countable, then we can simply take $g = f$. Hence, we may assume that $M$ is uncountable. Due to the Borel isomorphism theorem \cite[Theorem 15.6]{KechrisDST}, we can fix a Borel bijection $\R_{\geq 0} \to M \colon t \mapsto \phi_t$. For $t \in \R_{\geq 0} \cup \set{\infty}$, let $M_t \defeq \set{\phi_s \,:\, 0 \leq s < t}$. In particular, $M_0 = \0$ and $M_\infty = M$. Let $f_t \defeq f \recolor M_t$.
		
		Now let $\mu$ be a $G$-invariant probability measure. Observe the following facts:
		\begin{itemize}
			\item For each $\alpha \in \col{D}^+$, the function $t \mapsto \mu(f_t^{-1}(\alpha))$ is nondecreasing.
			
			\item For each $\beta \in \col{D}^-$, the function $t \mapsto \mu(f_t^{-1}(\beta))$ is nonincreasing.
			
			\item For each $\gamma \in \col{C} \setminus (\col{D}^+ \cup \col{D}^-)$, the function $t \mapsto \mu(f_t^{-1}(\gamma))$ is constant.
		\end{itemize}
		Furthermore, since $\mu$ is atomless, we have:
		\begin{itemize}
			\item For each $\gamma \in \col{C}$, the function $t \mapsto \mu(f^{-1}_t(\gamma))$ is continuous.
		\end{itemize}
		By definition, $f_\ast(\mu) \trianglelefteq (f_t)_\ast(\mu)$ if and only if either $f_\ast(\mu) = (f_t)_\ast(\mu)$, or there is $\alpha \in \col{D}^+$ such that
		\[
			\mu\left(f_t^{-1}(\alpha)\right) \,\leq\, \mu\left(f_t^{-1}(\beta)\right) \text{ for all } \beta \in \col{D}^-.
		\]
		This shows that the set of all $t \in \R_{\geq 0} \cup \set{\infty}$ such that $f_\ast(\mu) \trianglelefteq (f_t)_\ast(\mu)$ is closed. Hence, if we let \[
			T(\mu) \defeq \sup \set{t \in \R_{\geq 0} \cup \set{\infty} \,:\, f_\ast(\mu) \trianglelefteq (f_t)_\ast(\mu)} \qquad \text{and} \qquad g_\mu \defeq f_{T(\mu)},
		\]
		then $f_\ast(\mu) \trianglelefteq (g_\mu)_\ast(\mu)$. Also, by Lemma~\ref{lemma:dist2}, statements \ref{item:dist_norm} and \ref{item:norm_step} hold with $g_\mu$ in place of $g$.
		
		\begin{claim}
			The coloring $g_\mu$ is $(\mu, M)$-perfect.
		\end{claim}
		\begin{claimproof}
			We will show that, in fact, there is no recoloring move $\phi \in M$ that improves $g_\mu$. Indeed, suppose that $\phi \in M$ improves $g_\mu$. Then $\phi \not \in M_{T(\mu)}$, as otherwise we would have $\phi \subseteq g_\mu$. In particular, this means that $T(\mu) < \infty$.  Since  $\phi$ improves $g_\mu$, we can fix $\alpha \in \col{D}^+$ such that \[\mu\left(g_\mu^{-1}(\alpha)\right) \,<\, \mu\left(g_\mu^{-1}(\beta)\right) \text{ for all } \beta \in \col{D}^-.\] The functions $t \mapsto \mu(f_t^{-1}(\gamma))$, $\gamma \in \col{C}$, are continuous, so there is some $\epsilon > 0$ satisfying
		\[
		\mu\left(f_{T(\mu) + \epsilon}^{-1}(\alpha)\right) \,<\, \mu\left(f_{T(\mu) + \epsilon}^{-1}(\beta)\right)\text{ for all } \beta \in \col{D}^-.
		\]
		 This is a contradiction with the choice of $T(\mu)$.
		\end{claimproof}
	
		The above observations show that, with respect to the measure $\mu$, statements \ref{item:first}--\ref{item:last} hold with $g_\mu$ in place of $g$, which yields part \ref{item:step_a} of the lemma. To prove part \ref{item:step_b}, 
		we use the \hyperref[theo:UED]{Uniform Ergodic Decomposition Theorem}.
		 The statement is vacuous if $\Erg(G) = \0$, so assume that $\Erg(G) \neq \0$ and let $V \to \Erg(G) \colon x \mapsto \mu_x$ be a $G$-invariant Borel surjection given by Theorem~\ref{theo:UED}. Notice that, due to \cite[Theorem 17.25]{KechrisDST}, the function \[\Erg(G) \to \R_{\geq 0}\cup \set{\infty} \colon \mu \to T(\mu)\] is Borel. Hence, if we set $T(x) \defeq T(\mu_x)$ for each $x \in V$ and define \[g \colon V \to \col{C} \colon x \mapsto f_{T(x)}(x),\] then $g$ is a Borel $\col{C}$-coloring of $G$, and this coloring $g$ has all the desired properties. Indeed, for each $\mu \in \Erg(G)$, $g$ agrees with $g_\mu$ on the $\mu$-conull $G$-invariant set $\set{x \in V \,:\, \mu_x = \mu}$, and thus inherits the properties pertaining to $\mu$ from $g_\mu$.
	\end{scproof}

	\subsection{$\mu$-Equitability}\label{subsec:iteration}
	
	\begin{assum*}[for \S\ref{subsec:iteration}]
		Fix an aperiodic Borel graph $G$ of finite maximum degree $\Delta$ with vertex set $V$ and edge set $E$ and a finite set of colors $\col{C}$ of size $k \geq \Delta + 1$.
	\end{assum*}

	
	\begin{lemma}\label{lemma:all_measures}
		Let $f$ be a Borel proper $\col{C}$-coloring of $G$.
		\begin{enumerate}[label=\ep{\itshape\alph*}]
		    \item\label{item:all_measures_a} For every $G$-invariant probability measure $\mu$, there is a $\mu$\-/equitable $\col{C}$-coloring $g$ such that
		    \begin{equation}\label{eq:stable1}
		        \dist_\mu(f,g) \,\leq\, 7^{k+1} \cdot \disc_\mu(f).
		    \end{equation}
		    \item\label{item:all_measures_b} Furthermore, $G$ has a Borel proper $\col{C}$-coloring $g$ such that for every ergodic $G$-invariant probability measure $\mu$, $g$ is $\mu$-equitable and satisfies \eqref{eq:stable1}.
		\end{enumerate}
	\end{lemma}
	
	Note that, in view of the \hyperref[theo:UED]{Uniform Ergodic Decomposition Theorem} \ep{specifically, Theorem \ref{theo:UED}\ref{item:E2}}, the coloring $g$ given by Lemma~\ref{lemma:all_measures}\ref{item:all_measures_b} is in fact $\mu$-equitable for all $\mu \in \Inv(G)$. 
	
	\begin{scproof}\stepcounter{ForClaims} \renewcommand{\theForClaims}{\ref{lemma:all_measures}}
		The proofs of parts \ref{item:all_measures_a} and \ref{item:all_measures_b} of the lemma are virtually the same, except that the former relies on Lemma~\ref{lemma:step}\ref{item:step_a} and the latter on Lemma~\ref{lemma:step}\ref{item:step_b}. Below we give the proof of part \ref{item:all_measures_a} and highlight the only place where it needs to be modified to prove part \ref{item:all_measures_b}.
		
		We start by partitioning $\RM$ into countably many Borel $G$-separated sets.
		
		\begin{claim}\label{claim:color_moves}
			There is a Borel function $c \colon \RM \to \N$ such that for each $r \in \N$, $c^{-1}(r)$ is $G$-separated.
		\end{claim}
		\begin{claimproof}
			It will be more convenient to construct a Borel function $c \colon \RM \to \N^2$ such that for each pair $(r_0,r_1) \in \N^2$, $c^{-1}(r_0,r_1)$ is $G$-separated; this is sufficient as the set $\N^2$ is countable. With a slight \ep{but standard} abuse of notation, let $\fins{G}$ denote the set of all $S \in \fins{V}\setminus \set{\0}$ such that every two elements of $S$ are joined by a path in $G$. The proof of \cite[Lemma 7.3]{KechrisMiller} shows that there exists a Borel function $\theta \colon \fins{G} \to \N$ such that for each $r \in \N$, the sets in $\theta^{-1}(r)$ are pairwise disjoint. The mapping $c_0 \colon \RM \to \N \colon \phi \mapsto \theta(\dom(\phi) \cup N_G(\dom(\phi)))$ almost does the trick; the only problem is that two distinct recoloring moves $\phi$, $\psi$ may satisfy $\dom(\phi) \cup N_G(\dom(\phi)) = \dom(\psi) \cup N_G(\dom(\psi))$, in which case $c_0(\phi) = c_0(\psi)$. However, for each $S \in \fins{G}$, there are only finitely many $\phi \in \RM$ with $\dom(\phi) \cup N_G(\phi) = S$. Hence, by the Luzin--Novikov theorem \cite[Theorem 18.10]{KechrisDST}, there is a Borel function $c_1 \colon \RM \to \N$ such that if $\dom(\phi) \cup N_G(\dom(\phi)) = \dom(\psi) \cup N_G(\dom(\psi))$ for distinct $\phi$, $\psi \in \RM$, then $c_1(\phi) \neq c_1(\psi)$. Then we can let $c(\phi) \defeq (c_0(\phi), c_1(\phi))$.
		\end{claimproof}
	
		Now let $\mu$ be a $G$-invariant probability measure. We recursively construct a sequence of Borel proper $\col{C}$-colorings $f_n$, $n \in \N$, as follows. Fix a Borel function $c \colon \RM \to \N$ as in Claim \ref{claim:color_moves}. Let $(r_n, \col{D}^+_n, \col{D}^-_n)_{n \in \N}$ be a sequence of triples such that:
		\begin{enumerate}[label=\ep{\normalfont{}R\arabic*}]
			\item\label{item:type} for all $n\in \N$, $r_n \in \N$ and $\col{D}^+_n$, $\col{D}^-_n$ are disjoint nonempty subsets of $\col{C}$; and
			
			\item every triple $(r, \col{D}^+, \col{D}^-)$ as in \ref{item:type} 
			appears in the sequence $(r_n, \col{D}^+_n, \col{D}^-_n)_{n \in \N}$ infinitely often.
		\end{enumerate}
		Set $f_0 \defeq f$. Once $f_n$ is defined, let
		\[
			M_n \,\defeq\, \RM_3(f_n; \col{D}^+_n, \col{D}^-_n) \cap c^{-1}(r_n).
		\]
		Intersecting with $c^{-1}(r_n)$ ensures that $M_n$ is $G$-separated, so we can apply Lemma~\ref{lemma:step}\ref{item:step_a} with \[f_n,\ 3,\ \col{D}^+_n,\ \col{D}^-_n,\ \text{and}\ M_n \quad \text{in place of} \quad f,\ m,\ \col{D}^+,\ \col{D}^-,\ \text{and}\ M,\] in order to  obtain a Borel proper $\col{C}$-coloring $f_{n+1}$ such that: 
		\begin{enumerate}[label=\ep{\normalfont{}P\arabic*}]
			\item\label{item:dist_n} $\dist_\mu(f_n, f_{n+1}) \leq 3 \cdot \|(f_n)_\ast(\mu)-(f_{n+1})_\ast(\mu)\|_1$;
			\item\label{item:step_bound} for all $\alpha \in \col{D}^+((f_n)_\ast(\mu), (f_{n+1})_\ast(\mu))$, we have \[\|(f_n)_\ast(\mu)\,-\,(f_{n+1})_\ast(\mu)\|_1 \,\leq\, 6 \cdot (\mu(f_{n+1}^{-1}(\alpha)) \,-\, \mu(f_n^{-1}(\alpha)));\]
			\item\label{item:improve} $(f_n)_\ast(\mu) \trianglelefteq (f_{n+1})_\ast(\mu)$; and
			\item\label{item:perfect} $f_{n+1}$ is $(\mu, M_n)$-perfect.
		\end{enumerate}
		To establish part \ref{item:all_measures_b} of the lemma, we instead apply Lemma~\ref{lemma:step}\ref{item:step_b} to obtain a Borel proper $\col{C}$-coloring $f_{n+1}$ that satisfies \ref{item:dist_n}--\ref{item:perfect} for every $\mu \in \Erg(G)$. When the sequence $f_n$, $n \in \N$ is constructed, we define a Borel partial $\col{C}$-coloring $f_\infty \colon V \rightharpoonup \col{C}$ via the pointwise limit
		\[
		f_\infty(x) \,\defeq\, \lim_{n \to \infty} f_n(x).
		\]
		Since each $f_n$ is a proper coloring, $f_\infty$ is also proper, and since $k \geq \Delta + 1$, we can extend $f_\infty$ to a Borel proper $\col{C}$-coloring $g$ using Corollary~\ref{corl:extension}. We claim that this $g$ is as desired.
		
		To begin with, notice that conditions \ref{item:step_bound} and \ref{item:improve} enable us to apply Lemma~\ref{lemma:limit} to the sequence $((f_n)_\ast(\mu))_{n \in \N}$ with $A = 6$ and conclude, using \ref{item:dist_n}, that 
		\[
			\sum_{n=0}^\infty \dist_\mu(f_n, f_{n+1}) \,\leq\, 3 \cdot \sum_{n=0}^\infty \|(f_n)_\ast(\mu) \,-\, (f_{n+1})_\ast(\mu)\|_1 \,\leq\, \frac{7^{k+1}}{2} \cdot \disc_\mu(f) \,< \, \infty.
		\]
		By the Borel--Cantelli lemma, this implies that $f_\infty$ is defined for $\mu$-almost every $x \in V$; furthermore,
		\[
			\dist_\mu(f, g) \,=\, \dist_\mu(f, f_\infty) \,\leq\, \frac{7^{k+1}}{2} \cdot \disc_\mu(f).
		\]
		\ep{For simplicity, we bound the latter expression by $7^{k+1} \cdot \disc_\mu(f)$ in \eqref{eq:stable1}.} 
		It remains to show that $g$ is $\mu$-equitable. 
		To this end, we pass to a $\mu$-conull $G$-invariant Borel subset of $V$ and assume that $g = f_\infty$. Suppose, towards a contradiction,
		that $g$ is {not} $\mu$-equitable. Then, by Lemma~\ref{lemma:find_move}, $g$ is not $(\mu, \RM_3)$-perfect, i.e.,
		\[
			\mu\left(\set{\phi \in \RM_3 \,:\, \text{$\phi$ improves $g$ with respect to $\mu$}}\right) \,>\, 0.
		\]
		As there are only countably many triples $(r, \col{D}^+, \col{D}^-)$ with $r \in \N$ and $\col{D}^+$, $\col{D}^- \subseteq \col{C}$ disjoint and nonempty, we can choose $(r, \col{D}^+, \col{D}^-)$ so that
		\begin{equation}\label{eq:not_perfect}
		\mu\left(\set{\phi \in \RM_3(g; \col{D}^+, \col{D}^-) \cap c^{-1}(r) \,:\, \text{$\phi$ improves $g$ with respect to $\mu$}}\right) \,>\, 0.
		\end{equation}
		
		\begin{claim}\label{claim:pass}
			If $\phi \in \RM_3(g; \col{D}^+, \col{D}^-)$, then, for all sufficiently large $n \in \N$:
			\begin{enumerate}[label=\ep{\itshape\alph*}]
				\item\label{item:a} $\phi \in \RM_3(f_n; \col{D}^+, \col{D}^-)$;
				
				\item\label{item:b} if $\phi$ improves $g$ with respect to $\mu$, then $\phi$ also improves $f_n$ with respect to $\mu$.
			\end{enumerate}
		\end{claim}
		\begin{claimproof}
			\ref{item:a} This statement is implied by the fact that for all sufficiently large $n \in \N$, $f_n$ and $g$ agree on $\dom(\phi) \cup N_G(\dom(\phi))$ \ep{since, by our assumption, $g = f_\infty$}. 
			
			\ref{item:b} This follows from \ref{item:a} and the observation that if $\alpha$, $\beta \in \col{C}$ are such that $\mu(g^{-1}(\alpha)) < \mu(g^{-1}(\beta))$, then $\mu(f_n^{-1}(\alpha)) < \mu(f_n^{-1}(\beta))$ for all large enough $n \in \N$.
		\end{claimproof}
	
		Claim~\ref{claim:pass} and \eqref{eq:not_perfect} together show that there is $n_0 \in \N$ such that for all $n \geq n_0$,
		\[
			\mu\left(\set{\phi \in \RM_3(f_n;\col{D}^+, \col{D}^-) \cap c^{-1}(r) \,:\, \text{$\phi$ improves $f_n$ with respect to $\mu$}}\right) \,>\, 0.
		\]
		This precisely means that for all $n \geq n_0$,
		\[
			\text{$f_n$ is {not} $(\mu, \RM_3(f_n;\col{D}^+, \col{D}^-) \cap c^{-1}(r))$-perfect}.
		\]
		But this is a contradiction as there is $n \geq n_0$ with $(r_n, \col{D}^+_n, \col{D}^-_n) = (r, \col{D}^+, \col{D}^-)$, so
		\[
			\RM_3(f_n;\col{D}^+, \col{D}^-) \cap c^{-1}(r) \,=\, M_n,
		\]
		and $f_{n+1}$ is $(\mu, M_n)$-perfect by \ref{item:perfect}.
	\end{scproof}
	
	\subsection{Compressible graphs}\label{subsec:compressible}
	
	In this subsection we consider the case when $G$ is a Borel graph without any $G$-invariant probability measures; this situation is complementary to the one in Lemma~\ref{lemma:all_measures}. We begin by assembling some basic facts about such graphs. Throughout \S\ref{subsec:compressible}, $G$ denotes a locally countable Borel graph.
	
	A \ep{not necessarily finite} measure $\nu$ on a subset $U \subseteq V(G)$ is \emphd{$G$-invariant} if $\nu(A) = \nu(B)$ whenever $A$, $B \subseteq U$ and $A \approx_G B$. Note that $G$-invariance for a measure $\nu$ on $U \subseteq V(G)$ is in general stronger than $G[U]$-invariance. 
	There is a useful combinatorial characterization, due to Nadkarni~\cite{Nad} and subsequently generalized by Becker and Kechris \cite[Chapter 4]{BK}, of Borel sets $U \subseteq V(G)$ that do not support $G$-invariant probability measures, which we shall state after a few definitions. The \emphd{$G$-saturation} of a subset $U \subseteq V(G)$, denoted by $[U]_G$, is the union of all connected components of $G$ that intersect $U$. Observe that if $A$, $B \subseteq V(G)$ are Borel subsets such that $A \approx_G B$, then $[A]_G = [B]_G$. Since $G$ is locally countable, the Luzin--Novikov theorem \cite[Theorem 18.10]{KechrisDST} implies that $G$-saturations of Borel sets are Borel. A Borel set $U \subseteq V(G)$ is
	\begin{itemize}
		\item \emphd{$G$-compressible} if there is a Borel subset $A \subseteq U$ with $[A]_G = [U]_G$ and $U \approx_G U \setminus A$;
		
		\item \emphd{$G$-paradoxical} if there is a Borel partition $U = U_0 \sqcup U_1$ with $U_0 \approx_G U_1 \approx_G U$.
	\end{itemize}
	The graph $G$ itself is called \emphd{compressible} if $V(G)$ is a $G$-compressible set. Note that if $G$ is compressible, then every \emph{$G$-invariant} Borel subset $U \subseteq V(G)$ is $G$-compressible.
	
	\begin{theo}[{ess.~Nadkarni \cite{Nad}}]\label{theo:Nadkarni}
		Let $G$ be a locally countable Borel graph and let $U \subseteq V(G)$ be a Borel subset. The following statements are equivalent:
		\begin{enumerate}[label=\ep{\normalfont{}\roman*}]
			\item\label{item:comp} $U$ is not $G$-compressible;
			
			\item\label{item:paradox} $U$ is not $G$-paradoxical;
			
			\item\label{item:prob} there is a $G$-invariant probability measure $\nu$ on $U$;
			
			\item\label{item:fin} there is a $G$-invariant measure $\mu$ on $V(G)$ with $0 < \mu(U) < \infty$.
		\end{enumerate}
	\end{theo}
	\begin{scproof}
		Equivalence \ref{item:comp} $\Longleftrightarrow$ \ref{item:paradox} is proven in \cite[Proposition 2.1]{DJK}, Nadkarni's theorem \cite[Theorem 4.3.1]{BK} gives \ref{item:comp} $\Longleftrightarrow$ \ref{item:prob}, while \ref{item:prob} $\Longleftrightarrow$ \ref{item:fin} follows by \cite[Proposition 3.2]{DJK}.
	\end{scproof}


	For compressible graphs $G$, the relation $\approx_G$ can be understood quite well; for instance, Chen \cite{Ronnie} showed that if $G$ is compressible, the set of all $\approx_G$-equivalence classes of Borel subsets of $V(G)$ forms a \emph{cardinal algebra} in the sense of Tarski \cite{Tarski}. We will make use of the following:

	\begin{prop}[{\cite[Proposition 2.2]{DJK}}]\label{prop:sat}
		Let $G$ be a locally countable Borel graph and let $U \subseteq V(G)$ be a Borel subset. If $U$ is $G$-compressible, then $U \approx_G [U]_G$.
	\end{prop}

	\begin{corl}\label{corl:comp_exp}
		Let $G$ be a compressible Borel graph of finite maximum degree $\Delta$. If $U \subseteq V(G)$ is a Borel subset such that $N_G(U) \approx_G V(G)$, then $U \approx_G V(G)$ as well.
	\end{corl}
	\begin{scproof}
		Suppose, towards a contradiction, that $U \not \approx_G V(G)$. Proposition~\ref{prop:sat} then shows that, since $V(G) = [N_G(U)]_G \subseteq [U]_G$, the set $U$ is not $G$-compressible. By Theorem~\ref{theo:Nadkarni}, there is a $G$-invariant measure $\mu$ on $V(G)$ with $0 < \mu(U) < \infty$. But then \[0 \,<\, \mu(U) \,\leq\, \mu(V(G)) \,=\, \mu(N_G(U)) \,\leq\, \Delta \cdot \mu(U) \,<\, \infty,\] contradicting the compressibility of $G$.
	\end{scproof}

	\begin{corl}\label{corl:BEq}
		Let $G$ be a compressible locally countable Borel graph and let $f$ be a Borel proper $k$-coloring of $G$. The following statements are equivalent:
		\begin{enumerate}[label=\ep{\normalfont{}\roman*}]
			\item\label{item:eq} $f$ is Borel-equitable;
			
			\item\label{item:large} for every color $\alpha$, $f^{-1}(\alpha) \approx_G V(G)$;
			
			\item\label{item:full_comp} for every color $\alpha$, $[f^{-1}(\alpha)]_G = V(G)$ and $f^{-1}(\alpha)$ is $G$-compressible.
		\end{enumerate}
	\end{corl}
	\begin{scproof}
		Implication \ref{item:large} $\Longrightarrow$ \ref{item:eq} is clear, while \ref{item:full_comp} $\Longrightarrow$ \ref{item:large} follows by Proposition~\ref{prop:sat}. To prove \ref{item:eq} $\Longrightarrow$ \ref{item:full_comp}, let $f$ be Borel-equitable. This clearly implies that $[f^{-1}(\alpha)]_G = V(G)$ for every color $\alpha$, so it remains to show that $f^{-1}(\alpha)$ is $G$-compressible. Take any nonzero $G$-invariant measure $\mu$ on $V(G)$. Since $f$ is Borel-equitable, all the color classes of $f$ have the same $\mu$-measure, and hence
		\[
			\infty \,=\, \mu(V(G)) \,=\, k \cdot \mu(f^{-1}(\alpha)).
		\]
		This shows that $f^{-1}(\alpha)$ is $G$-compressible by Theorem~\ref{theo:Nadkarni}.
	\end{scproof}

	In addition to the equivalence relation $\approx_G$, it is useful to consider the preorder $\inj_G$, defined as follows: Given Borel sets $A$, $B \subseteq V(G)$, we write $A \inj_G B$ \ep{or, equivalently, $B \surj_G A$} if $A \approx_G B'$ for some Borel $B' \subseteq B$. Additionally, if $A \approx_G B'$ and for some Borel $B' \subseteq B$ such that $[B]_G = [B \setminus B']_G$, then we write $A \prec_G B$ \ep{or, equivalently, $B \succ_G A$}. Thus, in particular, a Borel set $U \subseteq V(G)$ is $G$-compressible if and only if $U \prec_G U$. Note that if $A \inj_G B$, then $[A]_G \subseteq [B]_G$ and $\mu(A) \leq \mu(B)$ for every $G$-invariant measure $\mu$. Furthermore, if $A \prec_G B$, then $\mu(A) < \mu(B)$ for every $G$-invariant measure $\mu$ such that $\mu([B]_G) > 0$. Indeed, if $A \approx_G B'$, where $B' \subseteq B$ satisfies $[B]_G = [B \setminus B']_G$, then $\mu(A) = \mu(B') = \mu(B) - \mu(B\setminus B')$, and if $\mu([B]_G) > 0$, then $\mu(B \setminus B') > 0$ as well.
	
	It is clear that the relation $\inj_G$ on Borel subsets of $V(G)$ is reflexive and transitive, and a standard Cantor--Schr\"oder--Bernstein-type argument shows that $A \approx_G B$ if and only if $A \inj_G B$ and $B \inj_G A$. While the relation $\inj_G$ generally fails to be a {total} preorder, any two Borel subsets of $V(G)$ can be made $\inj_G$-comparable by passing to suitable $G$-invariant subsets:
	
	\begin{prop}[{\cite[Lemma 4.5.1]{BK}}]\label{prop:compare}
		Let $G$ be a locally countable Borel graph and let $A$, $B \subseteq V(G)$ be Borel sets. Then there is a partition $V(G) = V_\prec \sqcup V_\approx \sqcup V_\succ$ of $V(G)$ into three $G$-invariant Borel subsets such that:
		\[
			A \cap V_\prec \prec_G B \cap V_\prec, \qquad A \cap V_\approx \approx_G B \cap V_\approx, \qquad \text{and} \qquad A \cap V_\succ \succ_G B \cap V_\succ.
		\]
	\end{prop}
	
	\begin{corl}\label{corl:meas_to_decomp}
		Let $G$ be a locally countable Borel graph and let $A$, $B \subseteq V(G)$ be Borel sets. If $\mu(A) = \mu(B)$ for all $\mu \in \Erg(G)$, then there is a partition $V(G) = V_0 \sqcup V_1$ of $V(G)$ into two $G$-invariant Borel sets such that $A \cap V_0 \approx_G B \cap V_0$ and $V_1$ is $G$-compressible.
	\end{corl}
	\begin{scproof}
		Let $V(G) = V_\prec \sqcup V_\approx \sqcup V_\succ$ be a partition given by Proposition~\ref{prop:compare} applied to $A$ and $B$. Consider any $\mu \in \Erg(G)$. Since $\mu$ is ergodic, precisely one of $V_\prec$, $V_\approx$, $V_\succ$ must be $\mu$-conull. Since $\mu(A) = \mu(B)$ by assumption, this implies that $\mu(V_\approx) = 1$ and $\mu(V_\prec \sqcup V_\succ) = 0$. As this holds for all $\mu \in \Erg(G)$, we conclude using Theorems~\ref{theo:UED} and~\ref{theo:Nadkarni} that the set $V_\prec \sqcup V_\succ$ is $G$-compressible. Hence, we can set $V_0 \defeq V_\approx$ and $V_1 \defeq V_\prec \sqcup V_\succ$.  
	\end{scproof}

	\begin{corl}\label{corl:compare}
		Let $G$ be a compressible locally countable Borel graph. Suppose that $U_0$, $U_1 \subseteq V(G)$ are Borel subsets with $U_0 \cup U_1 \approx_G V(G)$. Then there is a partition $V(G) = V_0 \sqcup V_1$ of $V(G)$ into two $G$-invariant Borel sets satisfying $U_0 \cap V_0 \approx_G V_0$ and $U_1 \cap V_1 \approx_G V_1$.
	\end{corl}
	\begin{scproof}
		Without loss of generality, we may assume that $U_0 \cup U_1 = V(G)$. From Proposition~\ref{prop:compare} we obtain a partition $V(G) = V_0 \sqcup V_1$ of $V(G)$ into two $G$-invariant Borel sets such that \[U_1 \cap V_0 \inj_G U_0 \cap V_0 \qquad \text{and}\qquad U_0 \cap V_1 \inj_G U_1 \cap V_1.\] We claim that these $V_0$, $V_1$ are as desired. Suppose that, say, $U_0 \cap V_0 \not \approx_G V_0$. Since $V(G) = U_0 \cup U_1$, the relation $U_1 \cap V_0 \inj_G U_0 \cap V_0$ implies $[U_0 \cap V_0]_G = V_0$. By Proposition~\ref{prop:sat}, the set $U_0 \cap V_0$ must be not $G$-compressible, so let $\mu$ be a $G$-invariant measure on $V(G)$ such that $0 < \mu(U_0 \cap V_0) < \infty$. Since $U_1 \cap V_0 \inj_G U_0 \cap V_0$, we conclude that
		\[
			0 \,<\, \mu(U_0 \cap V_0) \,\leq\, \mu(V_0) \,=\, \mu(U_0 \cap V_0) + \mu(U_1 \cap V_0) \,\leq\, 2\mu(U_0 \cap V_0) \,<\, \infty,
		\]
		contradicting the $G$-compressibility of $V_0$.
	\end{scproof}

	\begin{lemma}\label{lemma:two_ind}
		Let $G$ be a compressible Borel graph of finite maximum degree. Let $I$, $J \subseteq V(G)$ be disjoint Borel $G$-independent sets and suppose that $I \approx_G V(G)$. Then there exist disjoint Borel $G$-independent sets $I'$ and $J'$ such that $I' \cup J' = I \cup J$ and $I' \approx_G J' \approx_G V(G)$.
	\end{lemma}
	\begin{scproof}
		Applying Corollary~\ref{corl:compare} with $I \setminus N_G(J)$ and $I \cap N_G(J)$ in place of $U_0$ and $U_1$, we obtain a partition $V(G) = V_0 \sqcup V_1$ of $V(G)$ into two $G$-invariant Borel sets such that
		\[
			(I \setminus N_G(J)) \cap V_0 \approx_G V_0 \qquad \text{and} \qquad (I \cap N_G(J)) \cap V_1 \approx_G V_1.
		\]
		Since we can treat the induced subgraphs $G[V_0]$ and $G[V_1]$ separately, we may assume that either $V_0 = V(G)$ or $V_1 = V(G)$. Now we consider the two cases.
		
		If $V_0 = V(G)$ and $I \setminus N_G(J) \approx_G V(G)$, then, by Theorem~\ref{theo:Nadkarni}, the set $I \setminus N_G(J)$ is $G$-paradoxical, so there is a partition $I \setminus N_G(J) = I_0 \sqcup I_1$ of $I$ into two Borel sets satisfying $I_0 \approx_G I_1 \approx_G V(G)$. This allows us to set $I' \defeq I \setminus I_1 = (I \cap N_G(J)) \cup I_0$ and $J' \defeq J \cup I_1$.
			
		On the other hand, if $V_1 = V(G)$ and $I \cap N_G(J) \approx_G V(G)$, then $J \approx_G V(G)$ by Corollary~\ref{corl:comp_exp}, and hence we can simply take $I' \defeq I$ and $J' \defeq J$.
	\end{scproof}

	We are now ready to state and prove the main result of this subsection:

	\begin{theo}\label{theo:compressible}
		Let $G$ be a compressible Borel graph of finite maximum degree. If $\chi_\mathrm{B}(G) \leq k$, then $G$ has a Borel-equitable $k$-coloring.
	\end{theo}
	\begin{scproof}
		Let $\col{C}$ be a $k$-element set of colors and let $f \colon V(G) \to \col{C}$ be a Borel proper coloring of $G$. It follows from Corollary~\ref{corl:compare} that there is a partition $V(G) = \bigsqcup_{\alpha \in \col{C}} V_\alpha$ of $V(G)$ into $G$-invariant Borel sets satisfying $f^{-1}(\alpha) \cap V_\alpha \approx_G V_\alpha$ for each $\alpha \in \col{C}$. As we can deal with the induced subgraphs $G[V_\alpha]$ individually, we may assume that $V(G) = V_\alpha$ for some $\alpha \in \col{C}$ and thus $f^{-1}(\alpha) \approx_G V(G)$. Theorem~\ref{theo:compressible} then follows through a sequence of $k-1$ applications of Lemma~\ref{lemma:two_ind}.
	\end{scproof}

	We remark, incidentally, that the conclusion of Theorem~\ref{theo:compressible} may fail if the finite maximum degree assumption is replaced by local finiteness; the reason is that in a locally finite graph $G$, a set that is not $G$-compressible can still have a $G$-compressible neighborhood \ep{in contrast to Corollary~\ref{corl:comp_exp}}. We sketch a counterexample below. Take any aperiodic locally finite Borel graph $G$ such that $|\Erg(G)| = 1$ \ep{such graphs are called \emphd{uniquely ergodic}} and let $\mu$ be the unique ergodic $G$-invariant probability measure on $V(G)$. Then $\mu$ is atomless, so we can partition $V(G)$ as $V(G) = \bigsqcup_{n=1}^\infty V_n$, where each $V_n$ is a Borel set with $\mu(V_n) = 2^{-n}$. For each $x \in V(G)$, let $n(x) \in \N^+$ denote the index such that $x \in V_{n(x)}$, and let $H$ be the graph with vertex set \[V(H) \,\defeq\, \left\{(x, i) \,:\, x \in V(G) \text{ and } 1 \leq i \leq 2^{n(x)}\right\},\] in which two vertices $(x, i)$ and $(y,j)$ are adjacent if and only if $y \in \set{x} \cup N_G(x)$ and exactly one of $i$ and $j$ is equal to $1$. It is clear that $H$ is locally finite, and, owing to the fact that \[\int_{V(G)} 2^{n(x)} \,\Diff\mu(x) \,=\, \sum_{n=1}^\infty 2^n \cdot \mu(V_n) \,=\,  \infty,\] it is not hard to see that $H$ is compressible. Furthermore, $H$ has a Borel proper $2$-coloring, namely the function $V(H) \to \set{1,2} \colon (x,i) \mapsto \min\set{i, 2}$. Nevertheless, we claim that in every Borel proper $2$-coloring of $H$, one of the color classes must be not $H$-compressible \ep{and hence, by Corollary~\ref{corl:BEq}, such a coloring cannot be Borel-equitable}. Indeed, let $f \colon V(H) \to \set{1,2}$ be a Borel proper $2$-coloring of $H$. If $U \subseteq V(H)$ is a connected component of $H$, then $f$ must assign the same color to every vertex in $U \cap (V(G) \times \set{1})$ \ep{since any two such vertices are joined by a path of even length in $H$}. In other words, the function $V(G) \to \set{1,2} \colon x \mapsto f(x,1)$ is $G$-invariant. Since the measure $\mu$ is ergodic, there is a color $\alpha$ such that $f(x,1) = \alpha$ for a $\mu$-conull set of $x \in V(G)$. But this means that the pushforward of $\mu$ under the map $V(G) \to V(H) \colon x \mapsto (x,1)$ is an $H$-invariant probability measure on $f^{-1}(\alpha)$, showing that the color class $f^{-1}(\alpha)$ is not $H$-compressible.

	\subsection{Finishing the proof of Theorem~\ref{theo:stable_BHSz}}
	
	We are now in a position to complete the proof of Theorem~\ref{theo:stable_BHSz}. Part \ref{item:stable_a} of Theorem~\ref{theo:stable_BHSz} is given by Lemma~\ref{lemma:all_measures}\ref{item:all_measures_a}, so it remains to verify part \ref{item:stable_b}. To this end, we fix an aperiodic Borel graph $G$ of finite maximum degree $\Delta$ with vertex set $V$ and edge set $E$ and a finite set of colors $\col{C}$ of size $k \geq \Delta + 1$. Let $f \colon V \to \col{C}$ be a Borel proper coloring of $G$. By Lemma~\ref{lemma:all_measures}\ref{item:all_measures_b}, there exists a Borel proper coloring $h \colon V \to \col{C}$ such that for all $\mu \in \Erg(G)$, $h$ is $\mu$-equitable and
	\[\dist_\mu(f,h) \,\leq\, 7^{k+1} \cdot \disc_\mu(f).\]
	By applying Corollary~\ref{corl:meas_to_decomp} to each pair of color classes of $h$ and using the fact that finite \ep{and even countable} unions of $G$-compressible sets are $G$-compressible, we obtain a partition $V = V_0 \sqcup V_1$ of $V$ into two $G$-invariant Borel subsets such that:
	\begin{itemize}
		\item the sets $h^{-1}(\alpha) \cap V_0$, $\alpha \in \col{C}$ are pairwise $G$-equidecomposable; and
		
		\item the set $V_1$ is $G$-compressible.
	\end{itemize}
	By Theorem~\ref{theo:compressible}, the graph $G[V_1]$ has a Borel-equitable coloring $h' \colon V_1 \to \col{C}$. Define $g \colon V \to \col{C}$ by
	\[
		g(x) \,\defeq\, \begin{cases}
			h(x) &\text{if } x \in V_0;\\
			h'(x) &\text{if } x \in V_1.
		\end{cases}
	\]
	Then $g$ is a Borel-equitable $k$-coloring of $G$ such that for each $\mu \in \Erg(G)$,
	\[
		\dist_\mu(f,g) \,=\, \dist_\mu(f,h) \,\leq\, 7^{k+1} \cdot \disc_\mu(f),
	\]
	and the proof of Theorem~\ref{theo:stable_BHSz}\ref{item:stable_b} is complete.

	\section{Domination for partial colorings}\label{sec:domination}
	
	\subsection{Domination for list coloring}
	
	In this subsection, we prove Theorem~\ref{theo:list_domination}. Our strategy is to reduce Theorem~\ref{theo:list_domination}, through a series of auxiliary lemmas, to Theorem~\ref{theo:KN_domination}.
	
	
	\begin{lemma}\label{lemma:all_but_one}
		Let $G$ be a connected finite graph and let $\col{L}$ be a degree-list assignment for $G$. Suppose that $g$ is a partial proper $\col{L}$-coloring of $G$. Then for each vertex $u \in V(G)$, $G$ has a partial proper $\col{L}$-coloring $f$ with $\dom(f) \supseteq V(G) \setminus \set{u}$ and $f \dominates g$.
	\end{lemma}
	\begin{scproof}
		Suppose, towards a contradiction, that the tuple $(G, \col{L}, g, u)$ forms a counterexample that minimizes $|V(G)|$. Then clearly $|V(G)| \geq 2$. Without loss of generality, we may assume that the partial proper $\col{L}$-coloring $g$ is inclusion-maximal. This means that for each $x \in V(G) \setminus \dom(g)$, every color in $\col{L}(x)$ is assigned by $g$ to some neighbor of $x$. Since $|\col{L}(x)| \geq \deg_G(x)$, we conclude that $N_G(x) \subseteq \dom(g)$ and for each $\alpha \in \col{L}(x)$, there is exactly one $y \in N_G(x)$ with $g(y) = \alpha$.
		
		Consider now an arbitrary vertex $z \in V(G) \setminus \set{u}$ such that the subgraph $G - z$ is connected \ep{for example, if $T$ is a spanning subtree of $G$, then any leaf of $T$ distinct from $u$ is such}. If $z \not \in \dom(g)$, then pick any vertex $y \in N_G(z)$. Since $y$ is the unique neighbor of $z$ with the color $g(y)$, we may replace $g$ by the partial coloring $g^\ast$ with domain $(\dom(g) \setminus \set{y}) \cup \set{z}$  defined by
		\[
			g^\ast(x) \,\defeq\, \begin{cases}
				g(x) &\text{if } x \in \dom(g) \setminus \set{y};\\
				g(y) &\text{if } x = z.
			\end{cases}
		\]
		Therefore, we may assume that $z \in \dom(g)$. Let $g'$ be the restriction of $g$ onto $\dom(g) \cap (V(G) \setminus \set{z})$. For each $x \in V(G) \setminus \set{z}$, define
		\[
			\col{L}'(x) \,\defeq\, \begin{cases}
				\col{L}(x) \setminus \set{g(z)} &\text{if } x \in N_G(z);\\
				\col{L}(x) &\text{otherwise}.
			\end{cases}
		\]
		Then $\col{L}'$ is a degree-list assignment for $G-z$ and $g'$ is a partial proper $\col{L}'$-coloring. By the minimality of $|V(G)|$, $G - z$ has a partial proper $\col{L}'$-coloring $f'$ such that $\dom(f) \supseteq V(G) \setminus \set{z,u}$ and $f' \dominates g'$. But then the partial coloring $f \defeq f' \cup \set{(z, g(z))}$ satisfies the conclusion of the lemma.
	\end{scproof}

	\begin{lemma}\label{lemma:large_list}
		Let $G$ be a connected finite graph and let $\col{L}$ be a degree-list assignment for $G$. Suppose that $g$ is a partial proper $\col{L}$-coloring of $G$. If $G$ has no proper $\col{L}$-coloring $f$ with $f \dominates g$, then $|\col{L}(x)| = \deg_G(x)$ for all $x \in V(G)$.
	\end{lemma}
	\begin{scproof}
		Suppose that $x \in V(G)$ satisfies $|\col{L}(x)| \geq \deg_G(x) + 1$. Due to Lemma~\ref{lemma:all_but_one}, we may assume that $\dom(g) = V(G) \setminus \set{x}$. Then there is a color in $\col{L}(x)$ that is not assigned by $g$ to any of the neighbors of $x$, and hence $g$ can be extended to a proper $\col{L}$-coloring of $G$; a contradiction.
	\end{scproof}

	\begin{lemma}\label{lemma:same}
		Let $G$ be a connected finite graph without a cut-vertex and let $\col{L}$ be a degree-list assignment for $G$. Suppose that $g$ is a partial proper $\col{L}$-coloring of $G$. If $G$ has no proper $\col{L}$-coloring $f$ with $f \dominates g$, then all the lists $\col{L}(x)$, $x \in V(G)$, are equal to each other.
	\end{lemma}
	\begin{scproof}
		The statement is vacuous if $|V(G)| \leq 1$, so assume that $|V(G)| \geq 2$. Since $G$ is connected, it suffices to show that $\col{L}(x) = \col{L}(y)$ whenever $x$ and $y$ are adjacent. Suppose, towards a contradiction, that $x$ and $y$ are adjacent vertices and $\beta \in \col{L}(x) \setminus \col{L}(y)$. Due to Lemma~\ref{lemma:all_but_one}, we may assume that $\dom(g) = V(G) \setminus \set{x}$. Since $|\col{L}(x)| \geq \deg_G(x)$ and $g$ cannot be extended to a proper $\col{L}$-coloring of $G$, every color in $\col{L}(x)$ is assigned by $g$ to a single neighbor of $x$. In particular, there is a unique vertex $z \in N_G(x)$ with $g(z) = \beta$. Note that $z \neq y$ since $\beta \not \in \col{L}(y)$.
		
		Let $g'$ be the restriction of $g$ onto $\dom(g) \cap (V(G) \setminus \set{x, z})$. For each $u \in V(G) \setminus \set{x}$, define
		\[
			\col{L}'(u) \,\defeq\, \begin{cases}
				\col{L}(u) \setminus \set{\beta} &\text{if } u \in N_G(x);\\
				\col{L}(u) &\text{otherwise}.
			\end{cases}
		\]
		Then $\col{L}'$ is a degree-list assignment for $G-x$ and $g'$ is a partial proper $\col{L}'$-coloring. Furthermore, since $\beta \not \in \col{L}(y)$, we have $|\col{L}'(y)| \geq \deg_{G - x}(y) + 1$. Since $G$ has no cut-vertices, the graph $G - x$ is connected, so we may apply Lemma~\ref{lemma:large_list} to conclude that $G-x$ has a proper $\col{L}'$-coloring $f'$ such that $f' \dominates g'$. But then $f \defeq f' \cup \set{(x, \beta)}$ is a proper $\col{L}$-coloring of $G$ with $f \dominates g$; a contradiction.
	\end{scproof}

	We are ready to finish the proof of Theorem~\ref{theo:list_domination}. Let $G$ be a connected finite graph that is not a Gallai tree and let $\col{L}$ be a degree\-/list assignment for $G$. Let $g$ be a partial proper $\col{L}$-coloring of $G$ and suppose, towards a contradiction, that $G$ has no proper $\col{L}$-coloring $f$ with $f \dominates g$. Since $G$ is not a Gallai tree, $G$ has a block that is neither a clique nor an odd cycle. Let $U \subseteq V(G)$ be the vertex set of such a block and fix any $u \in U$. Due to Lemma~\ref{lemma:all_but_one}, we may assume that $\dom(g) = V(G) \setminus \set{u}$ and then replace $G$ by $G[U]$, $g$ by its restriction to $U$, and $\col{L}$ by the list assignment $\col{L}'(x) \defeq \col{L}(x) \setminus \set{g(y) \,:\, y \in N_G(x) \setminus U}$. In this way we arrange that $G$ is a graph without cut-vertices that is neither a clique nor an odd cycle. Let $\Delta$ be the maximum degree of $G$. By Lemma~\ref{lemma:same}, all the lists $\col{L}(x)$, $x \in V(G)$, are the same, and, by Lemma~\ref{lemma:large_list}, they all have size $\Delta$. Hence, if $\Delta \geq 3$, then we are done by Theorem~\ref{theo:KN_domination}. On the other hand, if $\Delta \leq 2$, then $G$ must be an even cycle. In that case, since $g$ is a proper $2$-coloring of the even path $G-u$, the neighbors of $u$ are colored the same by $g$, so $g$ can be extended to a proper $2$-coloring of $G$. 
	
	\subsection{One-ended subforests and measurable domination}
	
	Given a function $\phi$, we say that a sequence $x_0$, $x_1$, \ldots{} is \emphd{$\phi$-descending} if  $\phi(x_{n+1}) = x_n$ for all $n \in \N$. A function $\phi$ is \emphd{one-ended} if there is no infinite $\phi$-descending sequence. Let $G$ be a locally finite graph. Given a set $A \subseteq V(G)$, an \emphd{$A$-one-ended subforest} of $G$ is a one-ended function $\phi \colon V(G) \setminus A \to V(G)$ such that each $x \in V(G)\setminus A$ is adjacent to $\phi(x)$ in $G$. The word ``subforest'' is used here because if $\phi$ is an $A$-one-ended subforest of $G$, then the graph with vertex set $V(G)$ and edges joining each $x \in V(G) \setminus A$ to $\phi(x)$ is an acyclic subgraph of $G$. Given an $A$-one-ended subforest $\phi$ of $G$, 
	we define the \emphd{$\phi$-height} of a vertex $x \in V(G)$, in symbols $h_\phi(x)$, to be the greatest $n \in \N$ such that $x \in \im(\phi^n)$ \ep{such $n$ exists since $\phi$ is one-ended and $G$ is locally finite}. By definition, $h_\phi(x) = 0$ if and only if $x \not \in \im(\phi)$. Note that $h_\phi(\phi(x)) > h_\phi(x)$ for all $x \in V(G) \setminus A$.
	
	Conley, Marks, and Tucker-Drob developed the technique of one-ended subforests in order to prove a \hyperref[theo:meas_Brooks]{measurable version of Brooks's theorem} \cite{CMTD} \ep{see Theorem~\ref{theo:meas_Brooks}}. In particular, they showed that if $G$ is a Borel graph of finite maximum degree $\Delta$ and $A \subseteq V(G)$ is a Borel set such that $G$ has a Borel $A$-one-ended subforest, then $G$ has a Borel proper partial $\Delta$-coloring $f$ with $\dom(f) \supseteq V(G) \setminus A$ \cite[Lemma 3.9]{CMTD}. We strengthen this result by adding a domination requirement on $f$. Recall from \S\ref{subsec:compressible} that, given Borel sets $A$, $B \subseteq V(G)$, we write $A \inj_G B$ \ep{or, equivalently, $B \surj_G A$} if $A \approx_G B'$ for some Borel subset $B' \subseteq B$. For a pair of Borel partial colorings $f$, $g$ of $G$, we say that $f$ \emphd{Borel-dominates} $g$, in symbols $f \dominates_G g$, if $f^{-1}(\alpha) \surj_G g^{-1}(\alpha)$ for every color $\alpha$. Note that if $f \dominates_G g$, then $f \dominates_\mu g$ for every $G$-invariant probability measure $\mu$ on $V(G)$.
	
	\begin{lemma}\label{lemma:use_subf}
		Let $G$ be a Borel graph of finite maximum degree $\Delta$ and let $A \subseteq V(G)$ be a Borel set such that $G$ has a Borel $A$-one-ended subforest. If $g$ is a Borel proper partial $\Delta$-coloring of $G$, then $G$ has a Borel proper partial $\Delta$-coloring $f$ with $\dom(f) \supseteq V(G) \setminus A$ and $f \dominates_G g$.
	\end{lemma}
	\begin{scproof}\stepcounter{ForClaims} \renewcommand{\theForClaims}{\ref{lemma:use_subf}}
		Fix a Borel $A$-one-ended subforest $\phi$ of $G$. For each $n \in \N$, let \[V_n \defeq \set{x \in V(G) \,:\, h_\phi(x) = n} \qquad \text{and} \qquad V_{< n} \defeq \set{x \in V(G) \,:\, h_\phi(x) < n}.\] Let $\col{C}$ be a set of colors of size $\Delta$ and let $g$ be a Borel proper partial $\col{C}$-coloring of $G$. 
		We recursively construct a sequence of Borel proper partial $\col{C}$-colorings $(f_n)_{n=0}^\infty$, starting with $f_0 \defeq g$. We will ensure that each $f_n$ has the following properties:
		\begin{enumerate}[label=\ep{\normalfont{}D\arabic*}]
            \item\label{item:domain} $\dom(f_n) \supseteq V_{< n} \setminus A$;

            \item\label{item:stab} if $x \in V_{< n} \cap \dom(f_n)$, then $f_{n+1}(x) = f_{n}(x)$;	and
			
			\item\label{item:domination} for each color $\alpha \in \col{C}$ and every vertex $y \in f_n^{-1}(\alpha) \setminus f_{n+1}^{-1}(\alpha)$, there is $x \in f_{n+1}^{-1}(\alpha) \setminus (A \cup f_n^{-1}(\alpha))$ such that $h_\phi(x) = n$ and $\phi(x) = y$.
		\end{enumerate}
		Once $f_n$ is defined, we construct $f_{n+1}$ as follows. Let $f_n' \supseteq f_n$ be an arbitrary Borel inclusion-maximal proper partial $\col{C}$-coloring \ep{such $f_n'$ exists by Proposition~\ref{prop:list}}. By the maximality of $f_n'$, if $x \in V(G) \setminus \dom(f_n')$, then every neighbor of $x$ is colored by $f_n'$ and each color $\alpha \in \col{C}$ is used on precisely one neighbor of $x$. Hence, we may define $f_{n+1} \colon V(G) \rightharpoonup \col{C}$ by
		\[
		f_{n+1}(x) \,\defeq\, \begin{cases}
		f_n'(\phi(x)) &\text{if } x \in V_n \setminus (A \cup \dom(f_n'));\\
		\text{undefined} &\text{if } x \in \phi(V_n \setminus (A \cup \dom(f_n')));\\
		f_n'(x) &\text{otherwise}.
		\end{cases}
		\]
		Informally, we ``move'' the color from $\phi(x)$ to $x$ whenever $x \in V_n \setminus (A \cup \dom(f_n'))$. Conditions \ref{item:domain} and \ref{item:stab} are clearly satisfied. Notice also that the only vertices that are colored in $f_n$ but lose their color in $f_{n+1}$ are those of the form $\phi(x)$ for some $x \in V_n \setminus (A \cup \dom(f_n'))$, which implies \ref{item:domination}.
		
		The pointwise limit $f(x) \defeq \lim_{n \to \infty} f_n(x)$ is a proper partial $\col{C}$-coloring of $G$, and it follows from \ref{item:domain} and \ref{item:stab} that $\dom(f) \supseteq V(G) \setminus A$. It remains to show that $f \dominates_G g$. To this end, let
		\[
			X \defeq \set{x \in \dom(g) \,:\, f_n(x) = g(x) \text{ for all } n \in \N},
		\]
		and define a function $\psi \colon V(G) \to V(G)$ by setting
		\[
			\psi(x) \,\defeq\, \begin{cases}
			x &\text{if } x \in X \cup A;\\
			\phi(x) &\text{otherwise}.
			\end{cases}
		\]
		We claim that $\psi(f^{-1}(\alpha)) \supseteq g^{-1}(\alpha)$ for every $\alpha \in \col{C}$, which implies that $f \dominates_G g$ since the function $\psi$ has a Borel right inverse by the Luzin--Novikov theorem \cite[Theorem 18.10]{KechrisDST}. Let $\alpha \in \col{C}$ and consider any $y \in g^{-1}(\alpha)$. If $y \in X$, then $y \in f^{-1}(\alpha)$ and $y = \psi(y)$. Otherwise, there is $n \in \N$ such that $y \in f_n^{-1}(\alpha) \setminus f_{n+1}^{-1}(\alpha)$. Then, by~\ref{item:domination}, there is $x \in f_{n+1}^{-1}(\alpha) \setminus (A \cup f_n^{-1}(\alpha))$ with $h_\phi(x) = n$ and $\phi(x) = y$. Since $f_{n+1}(x) = \alpha$ but $x \not \in f_n^{-1}(\alpha)$, we conclude that $x \not \in X$, so $y = \phi(x) = \psi(x)$. Since $h_\phi(x) = n$, \ref{item:stab} yields $x \in f^{-1}(\alpha)$, and we are done. 
	\end{scproof}

	Given a Borel graph $G$ and a Borel subset $A \subseteq V(G)$, does $G$ have a Borel $A$-one-ended subforest? One case in which the answer is positive is when $A$ intersects every connected component of $G$:

	\begin{theo}[{Conley--Marks--Tucker-Drob \cite[Proposition 3.1]{CMTD}}]\label{theo:subf_full}
		Let $G$ be a locally finite Borel graph and suppose that $A \subseteq V(G)$ is a Borel subset that intersects every connected component of $G$. Then $G$ has a Borel $A$-one-ended subforest.
	\end{theo}

	Combining this with Lemma~\ref{lemma:use_subf} and Theorem~\ref{theo:list_domination}, we obtain the following:

	\begin{corl}\label{corl:Borel_domination}
		Let $G$ be a Borel graph of finite maximum degree $\Delta$. Suppose that least one of the following statements holds:
		\begin{enumerate}[label=\ep{\itshape\alph*}]
			\item\label{item:deg} every connected component of $G$ contains a vertex of degree less than $\Delta$; or
			
			\item\label{item:not_Gallai} no connected component of $G$ is a Gallai tree.
		\end{enumerate}
		If $g$ is a Borel proper partial $\Delta$-coloring of $G$, then $G$ has a Borel proper $\Delta$-coloring $f$ with $f \dominates_G g$.
	\end{corl}
	\begin{scproof}\stepcounter{ForClaims} \renewcommand{\theForClaims}{\ref{corl:Borel_domination}}
		\ref{item:deg} Let $A$ be the set of all $x\in V(G)$ with $\deg_G(x) < \Delta$. By assumption, $A$ intersects every connected component of $G$, so, by Theorem~\ref{theo:subf_full}, $G$ has a Borel $A$-one-ended subforest. Thus, due to Lemma~\ref{lemma:use_subf}, we may assume that $\dom(g) \supseteq V(G) \setminus A$. Let $f \supseteq g$ be a Borel inclusion-maximal proper partial $\Delta$-coloring \ep{such $f$ exists by Proposition~\ref{prop:list}}. The maximality of $f$ and the definition of $A$ imply that $A \subseteq \dom(f)$. Therefore, $\dom(f) = V(G)$, as desired.
		
		\ref{item:not_Gallai} This argument is essentially the same as the proof of \cite[Theorem 4.1]{CMTD}, with Lemma~\ref{lemma:use_subf} and Theorem~\ref{theo:list_domination} replacing \cite[Lemma 3.9]{CMTD} and Theorem~\ref{theo:list_Brooks} respectively. With a slight \ep{but standard} abuse of notation, let $\fins{G}$ denote the set of all $S \in \fins{V} \setminus \set{\0}$ such that every two elements of $S$ are joined by a path in $G$. Let $W \subseteq \fins{G}$ be the set of all $S \in \fins{G}$ such that $G[S]$ is connected and {not} a Gallai tree. Let $H$ be the graph with vertex set $W$ in which two distinct vertices $S$, $T$ are adjacent if and only if $(S \cup N_G(S)) \cap (T \cup N_G(T)) \neq \0$.
		
		\begin{claim}\label{claim:count_col}
			$\chi_\mathrm{B}(H) \leq \aleph_0$.
		\end{claim}
		\begin{claimproof}
			The proof of \cite[Lemma 7.3]{KechrisMiller} shows that there exists a Borel function $\theta \colon \fins{G} \to \N$ such that for each $r \in \N$, the sets in $\theta^{-1}(r)$ are pairwise disjoint. The map $c_0 \colon W \to \N \colon S \mapsto \theta(S \cup N_G(S))$ is almost a proper $\N$-coloring of $H$; the only problem is that two distinct sets $S$, $T \in W$ may satisfy $S \cup N_G(S) = T \cup N_G(T)$, in which case $c_0(S) = c_0(T)$. However, for each $F \in \fins{G}$, there are only finitely many $S \in W$ with $S \cup N_G(S) = F$. Hence, by the Luzin--Novikov theorem \cite[Theorem 18.10]{KechrisDST}, there is a Borel map $c_1 \colon W \to \N$ such that if $S \cup N_G(S) = T \cup N_G(T)$ for distinct $S$, $T \in W$, then $c_1(S) \neq c_1(T)$. Then $S \mapsto (c_0(S), c_1(S))$ is a Borel proper $\N^2$-coloring of $H$.
		\end{claimproof}
	
		From Claim~\ref{claim:count_col} and Corollary~\ref{corl:max_indep} we conclude that there is a Borel maximal $H$-independent set $I \subseteq W$. Let $A \defeq \bigcup I$. A connected graph is a Gallai tree if and only if all its finite connected subgraphs are Gallai trees; therefore, for every connected component $C$ of $G$, there is $S \in W$ such that $S \subseteq C$. Hence, the maximality of $I$ implies that $A$ intersects every connected component of $G$. By Theorem~\ref{theo:subf_full}, $G$ has a Borel $A$-one-ended subforest and thus, due to Lemma~\ref{lemma:use_subf},  we may assume that $\dom(g) \supseteq V(G) \setminus A$. For each vertex $x \in A$, let \[\col{L}(x) \defeq \col{C} \setminus \set{g(y) \,:\, y \in N_G(x) \cap (V(G) \setminus A) \cap \dom(g)}.\] Then $\col{L}$ is a degree-list assignment for the induced subgraph $G[A]$. For $S \in I$, let $g_S \colon S \rightharpoonup \col{C}$ denote the restriction of $g$ to $S$, so $g_S$ is a proper partial $\col{L}$-coloring of $G[S]$. Since $G[S]$ is a connected finite graph that is not a Gallai tree, by Theorem~\ref{theo:list_domination}, $G[S]$ admits a proper $\col{L}$-coloring $f_S \colon S \to \col{C}$ with $f_S \dominates g_S$. Furthermore, since there are only finitely many candidates for such $f_S$, the mapping $S \mapsto f_S$ can be arranged to be Borel. Now we can define a Borel proper coloring $f$ with $f \dominates_G g$ by
		\[
			f(x) \,\defeq\, \begin{cases}
				g(x) &\text{if } x \in V(G) \setminus A;\\
				f_S(x) &\text{if } x \in S \in I.
			\end{cases}
		\]
		\ep{Since $I$ is $H$-independent, the coloring $f$ is indeed well-defined and proper.}
	\end{scproof}

	To deal with graphs whose components are Gallai trees, we need another result of Conley, Marks, and Tucker-Drob. For a graph $G$, let $\mathsf{ic}(G)$ denote the number of infinite connected components of $G$ if it is finite, and $\infty$ otherwise. The \emphd{number of ends} of a connected locally finite graph $G$ is
	\[
		\mathsf{ends}(G) \,\defeq\, \sup \set{\mathsf{ic}(G[V(G)\setminus S]) \,:\, S \in \fins{V(G)}}. 
	\]
	If $\mathsf{ends}(G) = k$, then we say that $G$ is \emphd{$k$-ended}. Note that $\mathsf{ends}(G) = 0$ if and only if $G$ is finite.

	\begin{theo}[{Conley--Marks--Tucker-Drob \cite[proof of Theorem 4.2]{CMTD}}]\label{theo:subf_Gallai}
		Let $G$ be a locally finite Borel graph and let $\mu$ be a probability measure on $V(G)$. Suppose that every connected component $C \subseteq V(G)$ of $G$ has the following properties:
		\begin{itemize}
			\item $G[C]$ is a Gallai tree; and
			\item $\mathsf{ends}(G[C]) \not \in \set{0,2}$.
		\end{itemize}
		Then there is a $\mu$-conull $G$-invariant Borel subset $U \subseteq V(G)$ such that the induced subgraph $G[U]$ has a Borel $\0$-one-ended subforest.
	\end{theo}
	
	We are now ready to prove Theorem~\ref{theo:meas_domination}. 
	
	
	\begin{scproof}[ of Theorem~\ref{theo:meas_domination}]
		After passing to a $\mu$-conull $G$-invariant Borel subset of $V(G)$, we may assume that the partial coloring $g$ is in fact Borel. \ep{Since the leftover $\mu$-null set does not affect measurable domination, we may color it simply using the usual Brooks's theorem \cite[Theorem 5.2.4]{Die00}.} Partition $V(G)$ into three $G$-invariant Borel subsets as $V(G) = V_0 \sqcup V_1 \sqcup V_2$, where:
		\begin{itemize}
			\item every connected component of $G[V_0]$ has a vertex of degree less than $\Delta$;
			\item every vertex of $G[V_1]$ has degree $\Delta$ and no component of $G[V_1]$ is a Gallai tree;
			\item every vertex of $G[V_2]$ has degree $\Delta$ and every component of $G[V_2]$ is a Gallai tree.
		\end{itemize}
		For each $i \in \set{0,1,2}$, let $g_i$ denote the restriction of $g_i$ to $V_i$. By Corollary~\ref{corl:Borel_domination}, there exist Borel proper $\Delta$-colorings $f_0$ and $f_1$ of $G[V_0]$ and $G[V_1]$ respectively such that $f_0 \dominates_G g_0$ and $f_1 \dominates_G g_1$.
		
		Now consider the graph $G[V_2]$. Recall that a locally finite graph is \emphd{regular} if all its vertices have the same degree \ep{so the graph $G[V_2]$ is regular}. Observe that the only regular $0$-ended Gallai trees are cliques and odd cycles, while the only regular $2$-ended Gallai trees are two-way infinite paths. This implies that, since $\Delta \geq 3$ and $G$ does not contain a clique on $\Delta + 1$ vertices, $\mathsf{ends}(G[C]) \not \in \set{0,2}$ for every connected component $C$ of $G[V_2]$. By Theorem~\ref{theo:subf_Gallai}, after discarding a $\mu$-null $G$-invariant Borel set, we may assume that $G[V_2]$ admits a Borel $\0$-one-ended subforest. Hence, by Lemma~\ref{lemma:use_subf}, there is a Borel proper $\Delta$-coloring $f_2$ of $G[V_2]$ with $f_2 \dominates_G g_2$. Then $f \defeq f_0 \cup f_1 \cup f_2$ is a Borel proper $\Delta$-coloring of $G$ with $f \dominates_G g$ \ep{and hence also $f \dominates_\mu g$}, and we are done.
	\end{scproof}
	
	\section{Equitable $\Delta$-colorings}\label{sec:meas_KN}
	
	\subsection{Preliminary lemmas}
	
	In this section we prove Theorem~\ref{theo:meas_KN}. Our argument is analogous to the proof of Theorem~\ref{theo:KN} in \cite{KN}, modulo the changes necessary to adapt it to the measurable setting. In particular, the \hyperref[theo:HSz]{Hajnal--Szemer\'edi theorem} and Theorem~\ref{theo:KN_domination} are replaced by Theorems~\ref{theo:BHSz} and \ref{theo:meas_domination}, respectively. \ep{In fact, some of the final calculations presented in \S\ref{subsec:Delta} end up being somewhat simpler than the corresponding calculations in \cite{KN}.}
	
	We start by collecting a few preliminary results. First, we need a version of Theorem~\ref{theo:BHSz} for graphs that may have finite components:
	
	\begin{lemma}\label{lemma:meas_HSz}
		Let $G$ be a Borel graph of finite maximum degree $\Delta$ and let $\mu$ be an atomless $G$-invariant probability measure on $V(G)$. If $k \geq \Delta + 1$, then $G$ has a $\mu$-equitable $k$-coloring.
	\end{lemma}
	\begin{scproof}
		Let $U \subseteq V(G)$ be the union of all the infinite components of $G$. Then $U$ is a $G$-invariant Borel set, and, by Theorem~\ref{theo:BHSz}, $G[U]$ has a Borel-equitable $k$-coloring. If $\mu(U) = 1$, then we are done, so assume that $\mu(U) < 1$. Then, upon passing to the subgraph $G[V(G) \setminus U]$ and scaling $\mu$ appropriately, we may assume that every component of $G$ is finite.
		
		Let $\col{C}$ be a set of colors of size $k$. By Theorem~\ref{theo:KST}, $G$ has a Borel proper $\col{C}$-coloring $g$. For each function $\theta \colon \col{C} \to \N$, let $V_\theta \subseteq V(G)$ be the union of all the components $C$ of $G$ satisfying
		\[
			|\set{x \in C \,:\, g(x) = \alpha}| \,=\, \theta(\alpha) \qquad \text{for all } \alpha \in \col{C}.
		\]
		Then $V(G) = \bigsqcup \set{V_\theta \,:\, \theta \colon \col{C} \to \N}$ is a partition of $G$ into countably many $G$-invariant Borel sets. Again, whenever $\mu(V_\theta) \neq 0$, we may pass to the subgraph $G[V_\theta]$ and scale $\mu$ appropriately, thus reducing the situation to the case when $V(G) = V_\theta$ for some fixed $\theta \colon \col{C} \to \N$. Let $\mathsf{Sym}(\col{C})$ denote the set of all bijections $\col{C} \to \col{C}$. Since $\mu$ is atomless and every component of $G$ is finite, we can partition $V(G)$ into Borel $G$-invariant sets as $V(G) = \bigsqcup \set{V_\pi \,:\, \pi \in \mathsf{Sym}(\col{C})}$ so that $\mu(V_\pi) = 1/k!$ for all $\pi \in \mathsf{Sym}(\col{C})$. Then the coloring $f$ that sends each $x \in V_\pi$ to $(\pi \circ g)(x)$ is $\mu$-equitable.
	\end{scproof}

	To state our next lemma we need to introduce some terminology. Let $G$ be a Borel graph of finite maximum degree and let $\mu$ be a $G$-invariant probability measure on $V(G)$. Given a Borel subset $X \subseteq V(G)$, we define the \emphd{cost} of $X$ relative to $G$ and $\mu$ by the formula
	\begin{equation}\label{eq:cost}
		\mathsf{C}_\mu(G; X) \,\defeq\, \int_X |N_G(x) \setminus X| \,\Diff \mu(x) \,+\, \frac{1}{2}\int_X |N_G(x) \cap X| \,\Diff \mu(x).
	\end{equation}
	Intuitively, $\mathsf{C}_\mu(G;X)$ represents the ``normalized'' number of edges of $G$ incident to a vertex in $X$; the second summand in \eqref{eq:cost} is halved since the edges joining {two} vertices of $X$ are counted twice. In particular, we have $d_\mu(G) = 2\mathsf{C}_\mu(G; V(G))$ \ep{recall that $d_\mu(G)$ is the $\mu$-average degree of $G$}. The $G$-invariance of $\mu$ implies that if $X$, $Y \subseteq V(G)$ are disjoint Borel sets, then
	\begin{equation}\label{eq:cost_sum}
		\mathsf{C}_\mu(G;X \sqcup Y) \,=\, \mathsf{C}_\mu(G;X) \,+\, \mathsf{C}_\mu(G; Y) \,-\, \int_Y |N_G(y) \cap X| \,\Diff \mu(y).
	\end{equation}
	Note also that if $X \subseteq Y$, then $\mathsf{C}_\mu(G;X) \leq \mathsf{C}_\mu(G;Y)$.
	
	\begin{lemma}\label{lemma:X}
		Let $G$ be a Borel graph of finite maximum degree and let $\mu$ be a $G$\-/invariant probability measure on $V(G)$. Then, for each real number $t \geq 0$, there exists a Borel subset $X \subseteq V(G)$ with the following properties:
		\begin{enumerate}[label=\ep{\normalfont{}X\arabic*}]
			\item\label{item:X1} $\deg_G(y) < 2t$ for all $y \in V(G) \setminus X$;
			\item\label{item:X2} $|N_G(y) \setminus X| < t$ for all $y \in V(G) \setminus X$; and
			\item\label{item:X3} $\mathsf{C}_\mu(G; X') \geq t \mu(X')$ for every Borel set $X' \subseteq X$.
		\end{enumerate}
	\end{lemma}
	\begin{scproof}
		By Theorem~\ref{theo:KST}, $\chi_\mathrm{B}(G)$ is finite, so fix a Borel proper coloring $c \colon V(G) \to \set{0, \ldots, k-1}$ for some $k \in \N^+$. Recursively construct Borel sets $X_r$, $0 \leq r \leq k$, as follows: Set \[X_0 \defeq \set{x \in V(G) \,:\, \deg_G(x) \geq 2t},\] and, once $X_r$ is defined for some $0 \leq r < k$, let
		\[
			Y_r \defeq \set{y \in V(G) \setminus X_r \,:\, |N_G(y)\setminus X_r| \geq t}, \qquad I_r \defeq Y_r \cap c^{-1}(r), \qquad \text{and} \qquad X_{r+1} \defeq X_r \cup I_r.
		\]
		We claim that the set $X \defeq X_k$ is as desired. Condition~\ref{item:X1} is a consequence of the definition of $X_0$. To show \ref{item:X2}, suppose that $y \in V(G) \setminus X$ satisfies $|N_G(y) \setminus X| \geq t$ and let $r \defeq c(y)$. Then $|N_G(y) \setminus X_r| \geq |N_G(y) \setminus X| \geq t$ and hence $y \in Y_r$. But then $y \in Y_r \cap c^{-1}(r) = I_r \subseteq X$; a contradiction.
		It remains to verify \ref{item:X3}. Let $X' \subseteq X$ be a Borel subset. For each $0 \leq r < k$, let $X_r' \defeq X' \cap X_r$ and $I'_r \defeq X' \cap I_r$. By repeatedly applying \eqref{eq:cost_sum}, we obtain
		\[
			\mathsf{C}_\mu(G; X') \,=\, \mathsf{C}_\mu(G; X_0') \,+\, \sum_{r=0}^{k-1} \left(\mathsf{C}_\mu(G; I_r') \,-\, \int_{I_r'} |N_G(y) \cap X_r'| \,\Diff\mu(y) \right).
		\]
		From \eqref{eq:cost} and the definition of $X_0$, it follows that
		\[
			\mathsf{C}_\mu(G;X_0') \,\geq\, \frac{1}{2} \int_{X_0'} \deg_G(x) \,\Diff\mu(x) \,\geq\, t\mu(X_0').
		\]
		Let $0 \leq r < k$. The set $I_r'$ is $G$-independent, so
		\[
			\mathsf{C}_\mu(G; I_r') \,=\, \int_{I_r'} \deg_G(y) \,\Diff\mu(y).
		\]
		Therefore, we have
		\[
			\mathsf{C}_\mu(G; I_r') \,-\, \int_{I_r'} |N_G(y) \cap X_r'| \,\Diff\mu(y) \,=\, \int_{I_r'} |N_G(y) \setminus X_r'| \,\Diff\mu(y) \,\geq\, \int_{I_r'} |N_G(y) \setminus X_r| \,\Diff\mu(y) \,\geq\, t\mu(I_r').
		\]
		Putting everything together, we obtain the desired inequality
		\[
			\mathsf{C}_\mu(G; X') \,\geq\, t\mu(X_0') \,+\, \sum_{r=0}^{k-1} t \mu(I_r') \,=\, t\mu(X'). \qedhere
		\]
	\end{scproof}

	We need one more technical lemma:

	\begin{lemma}\label{lemma:quick_move}
		Let $G$ be a locally finite Borel graph and let $\mu$ be an atomless $G$\=/invariant probability measure on $V(G)$. Let $\col{C}$ be a finite set of colors and let $f$ be a Borel proper $\col{C}$-coloring of $G$. Then, for every Borel set $X \subseteq V(G)$, there is a Borel proper $\col{C}$-coloring $g$ of $G$ such that:
		\begin{enumerate}[label=\ep{\normalfont{}P\arabic*}]
			\item $g(x) = f(x)$ for all $x \in X$; and
			\item\label{item:P2} for all $\alpha$, $\beta \in \col{C}$, if $\mu(g^{-1}(\alpha)) < \mu(g^{-1}(\beta))$, then
			\[
				\mu(\set{y \in g^{-1}(\beta) \setminus X \,:\, N_G(y) \cap g^{-1}(\alpha) = \0}) \,=\, 0.
			\]
		\end{enumerate}
	\end{lemma}
	\begin{scproof}\stepcounter{ForClaims} \renewcommand{\theForClaims}{\ref{lemma:quick_move}}
		This is a \ep{significantly simpler} variant of the proof of Lemma \ref{lemma:all_measures}. Let $(r_n, \alpha_n, \beta_n)_{n \in \N}$ be a sequence of triples such that:
		\begin{enumerate}[label=\ep{\normalfont{}R\arabic*}]
			\item\label{item:type1} for all $n\in \N$, $r_n \in \N$ and $\alpha_n$, $\beta_n \in \col{C}$ are distinct colors; and
			
			\item every triple $(r, \alpha, \beta)$ as in \ref{item:type1} 
			appears in the sequence $(r_n, \alpha_n, \beta_n)_{n \in \N}$ infinitely often.
		\end{enumerate}
		Fix a Borel proper coloring $c \colon V(G) \to \N$ of $G$ \ep{for instance, we could take $c = f$}. Recursively construct Borel proper $\col{C}$-colorings $f_n$, $n \in \N$, as follows. Set $f_0 \defeq f$. Once $f_n$ is defined, we split the definition of $f_{n+1}$ into two cases.
		
		\begin{leftbar}
		\noindent \textbf{Case 1:} $\mu(f_n^{-1}(\alpha_n)) \geq \mu(f_n^{-1}(\beta_n))$. Then set $f_{n+1} \defeq f_n$.
		
		\medskip
		    
		\noindent \textbf{Case 2:} $\mu(f_n^{-1}(\alpha_n)) < \mu(f_n^{-1}(\beta_n))$. Let
		\[
			A_n \,\defeq\, \set{y \in f_n^{-1}(\beta_n) \setminus X \,:\, N_G(y) \cap f_n^{-1}(\alpha_n) = \0} \cap c^{-1}(r_n),
		\]
		
		\begin{leftbar}
		\noindent \textbf{Subcase 1.1:} $\mu(f_n^{-1}(\alpha_n)) + \mu(A_n) \leq \mu(f_n^{-1}(\beta_n)) - \mu(A_n)$. Then define $B_n \defeq A_n$.
		
		\medskip
		
		\noindent \textbf{Subcase 2.2:} $\mu(f_n^{-1}(\alpha_n)) + \mu(A_n) > \mu(f_n^{-1}(\beta_n)) - \mu(A_n)$. Since $\mu$ is atomless, we can then let $B_n \subseteq A_n$ be an arbitrary Borel subset of $A_n$ with $\mu(B_n) = (\mu(f_n^{-1}(\beta_n)) - \mu(f_n^{-1}(\alpha_n)))/2$.
		\end{leftbar}
		
		\noindent Note that in both subcases we have
		\begin{equation}\label{eq:BMeasure}
			\mu(B_n) \,=\, \min \left\{A_n,\, \frac{\mu(f_n^{-1}(\beta_n)) - \mu(f_n^{-1}(\alpha_n))}{2}\right\}.
		\end{equation}
		To finish Case 2, define
		\[
		f_{n+1}(x) \,\defeq\, \begin{cases}
		f_n(x) &\text{if } x \in V(G) \setminus B_n;\\
		\alpha_n &\text{if } x \in B_n.
		\end{cases}
		\]
		By construction, $f_{n+1}$ is proper and $f_{n+1}(x) = f_n(x) = f(x)$ for all $x \in X$.
		\end{leftbar}
		
		
		Now we define a Borel partial $\col{C}$-coloring $f_\infty \colon V \rightharpoonup \col{C}$ via the pointwise limit
		\[
		f_\infty(x) \,\defeq\, \lim_{n \to \infty} f_n(x).
		\]
		Since each $f_n$ is a proper coloring, $f_\infty$ is also proper. We wish to show that $f_\infty$ is defined $\mu$-almost everywhere. While this fact can be derived using Lemma~\ref{lemma:limit}, in this case we can give a simpler and more straightforward convergence argument. Let us introduce the following notation:
		\[
			\omega_n \colon \col{C} \to [0,1] \colon \gamma \mapsto \mu(f^{-1}_n(\gamma)) \qquad \text{and} \qquad S_n \,\defeq\, \frac{1}{2}\sum_{\gamma \in \col{C}} \sum_{\delta \in \col{C}}  \left|\omega_n(\gamma) - \omega_n(\delta)\right|.
		\]
		
		\begin{claim}\label{claim:semiinvariant}
			For all $n \in \N$, we have $\dist_\mu(f_n, f_{n+1}) \leq (S_n - S_{n+1})/2$.
		\end{claim}
		\begin{claimproof}
			If in the construction of $f_{n+1}$ Case 1 occurred, then $f_{n+1} = f_n$, and hence both sides of the desired inequality are $0$. Now assume that Case 2 occurred. Observe that if $a$, $b$, $c$, $d$ are real numbers with $0 \leq d \leq (b-a)/2$, then
			\begin{equation}\label{eq:ineq}
				\left|c - a\right| \,+\, \left|c - b\right| \,-\, \left|c - a - d\right| \,-\, \left|c - b + d\right| \,\geq\, 0.
			\end{equation}
			By construction, $\dist_\mu(f_n, f_{n+1}) = \mu(B_n)$. For each $\gamma \in \col{C}$, we have
			\[
				\omega_{n+1}(\gamma) \,=\, \begin{cases}
					\omega_n(\alpha_n) + \mu(B_n) &\text{if } \gamma = \alpha_n;\\
					\omega_n(\beta_n) - \mu(B_n) &\text{if } \gamma = \beta_n;\\
					\omega_n(\gamma) &\text{otherwise}.
				\end{cases}
			\]
			It follows from \eqref{eq:BMeasure} and \eqref{eq:ineq} that, for each $\gamma \in \col{C}\setminus \set{\alpha_n, \beta_n}$,
			\[
			\left|\omega_n(\gamma) - \omega_n(\alpha_n)\right| \,+\, \left|\omega_n(\gamma) - \omega_n(\beta_n)\right| \,-\, \left|\omega_{n+1}(\gamma) - \omega_{n+1}(\alpha_n)\right| \,-\, \left|\omega_{n+1}(\gamma) - \omega_{n+1}(\beta_n)\right| \,\geq\, 0.
			\]
			Therefore,
			\[
				S_n - S_{n+1} \,\geq\, \left|\omega_n(\alpha_n) - \omega_n(\beta_n)\right| \,-\, \left|\omega_{n+1}(\alpha_n) - \omega_{n+1}(\beta_n)\right| \,=\, 2\mu(B_n) \,=\, 2\dist_\mu(f_n, f_{n+1}),
			\]
			as desired.
		\end{claimproof}
	
		Claim~\ref{claim:semiinvariant} yields $\sum_{n=0}^\infty \dist_\mu(f_n, f_{n+1}) \leq S_0/2 < \infty$, and hence, by the Borel--Cantelli lemma, the domain of $f_\infty$ is $\mu$-conull. Thus, there is a $\mu$-conull $G$-invariant Borel subset $U \subseteq \dom(f_\infty)$. Define $g \colon V(G) \to \col{C}$ by sending each $x \in U$ to $f_\infty(x)$ and each $x \in V(G)\setminus U$ to $f(x)$. We claim that this $g$ is as desired. It is clear that $g$ is proper and that $g(x) = f(x)$ for all $x \in X$. 
		
		It remains to verify \ref{item:P2}. To this end, let $\alpha$, $\beta \in \col{C}$ be colors such that $\mu(g^{-1}(\alpha)) < \mu(g^{-1}(\beta))$. We shall argue that every vertex $y \in U \cap (g^{-1}(\beta) \setminus X)$ has a neighbor in $g^{-1}(\alpha)$, which implies \ref{item:P2} since $\mu(U)=1$. Suppose, towards a contradiction, that $y \in U \cap (g^{-1}(\beta) \setminus X)$ satisfies $N_G(y) \cap g^{-1}(\alpha) = \0$. Set $r \defeq c(y)$. Since $\set{y} \cup N_G(y) \subseteq U$, there is $n_0 \in \N$ such that for all $n \geq n_0$, we have:
		\begin{enumerate}[label=\ep{\normalfont{}L\arabic*}]
			\item\label{item:L1} $f_n(y) = \beta$ and $N_G(y) \cap f_n^{-1}(\alpha) = \0$; and
			\item\label{item:L2} $\mu(f_n^{-1}(\alpha)) < \mu(f_n^{-1}(\beta))$.
		\end{enumerate}
		Take any $n \geq n_0$ with $(r_n, \alpha_n, \beta_n) = (r, \alpha, \beta)$. It follows from \ref{item:L2} that in the construction of $f_{n+1}$, Case 2 occurred. Also, by \ref{item:L1}, we have $y \in A_n$. On the other hand, $y \not \in B_n$ since $f_{n+1}(y) = \beta \neq \alpha$. Therefore, $B_n \neq A_n$, which implies that Subcase 2.2 occurred in the construction of $f_{n+1}$. Hence,
		\[
			\mu(f_n^{-1}(\alpha)) + \mu(B_n) \,=\, \mu(f_n^{-1}(\beta)) - \mu(B_n),
		\]
		i.e., $\mu(f_{n+1}^{-1}(\alpha)) = \mu(f_{n+1}^{-1}(\beta))$. But this contradicts \ref{item:L2} with $n+1$ in place of $n$.
	\end{scproof}
	
	\subsection{Proof of Theorem~\ref{theo:meas_KN}}\label{subsec:Delta}
	
	We are now fully equipped to prove Theorem~\ref{theo:meas_KN}, so let $G$ be a Borel graph of finite maximum degree $\Delta \geq 3$ without a clique on $\Delta+1$ vertices, and let $\mu$ be an atomless $G$-invariant probability measure on $V(G)$ such that $d_\mu(G) \leq \Delta/5$. First we use Lemma~\ref{lemma:X} with $t = 2\Delta/5$ to find a Borel subset $X \subseteq V(G)$ such that:
	\begin{enumerate}[label=\ep{\normalfont{}X\arabic*}]
		\item\label{item:X11} $\deg_G(y) < 4\Delta/5$ for all $y \in V(G) \setminus X$;
		\item\label{item:X21} $|N_G(y) \setminus X| < 2\Delta/5$ for all $y \in V(G) \setminus X$; and
		\item\label{item:X31} $\mathsf{C}_\mu(G; X') \geq (2\Delta/5) \mu(X')$ for every Borel set $X' \subseteq X$.
	\end{enumerate}

	\begin{big_claim}\label{claim:X_small}
		$\mu(X) \leq 1/4$.
	\end{big_claim}
	\begin{scproof}
		This is a consequence of the following chain of inequalities:
		\[
			\frac{2\Delta \mu(X)}{5} \,\leq\, \mathsf{C}_\mu(G;X) \,\leq\,\mathsf{C}_\mu(G; V(G)) \,=\, \frac{d_\mu(G)}{2} \,\leq\, \frac{\Delta}{10}.\qedhere
		\]
	\end{scproof}

	Next we apply Lemma~\ref{lemma:meas_HSz} to obtain a $\mu$-equitable $(\Delta+1)$-coloring $h$ of the subgraph $G[X]$; here ``$\mu$-equitable'' means that each color class of $h$ has measure $\mu(X)/(\Delta+1)$. As in the proof of Corollary~\ref{corl:almost_equit}, we then uncolor one of the color classes of $h$ and apply Theorem~\ref{theo:meas_domination} to the resulting partial coloring. 
	This produces a $\mu$-measurable proper $\Delta$-coloring $h^\ast$ of $G[X]$ in which every color class has measure {at least} $\mu(X)/(\Delta+1)$. After passing to a $\mu$-conull $G$-invariant Borel subset of $V(G)$, we may assume that $h^\ast$ is Borel.
	
	\begin{big_claim}\label{claim:h_bounds}
	If $S \subseteq X$ is the union of some $s$ color classes of $h^\ast$, then
	\[
		\frac{s\mu(X)}{\Delta + 1} \,\leq\, \mu(S) \,\leq\, \frac{(s+1) \mu(X)}{\Delta + 1}.
	\]
	\end{big_claim}
	\begin{scproof}
		Immediate from the fact that each color class of $h^\ast$ has measure at least $\mu(X)/(\Delta + 1)$.
	\end{scproof}

	Now we use Proposition~\ref{prop:list} to obtain a Borel inclusion-maximal proper partial $\Delta$-coloring $g \supseteq h^\ast$ of $G$. Then $X \subseteq \dom(g)$ by definition; on the other hand, every vertex in $V(G) \setminus X$ has degree less than $4\Delta/5 < \Delta$, so $V(G) \setminus X \subseteq \dom(g)$ as well. In other words, $\dom(g) = V(G)$, i.e., $g$ is a {total} coloring. Finally, we invoke Lemma~\ref{lemma:quick_move} to produce a Borel proper $\Delta$-coloring $f$ of $G$ such that:
	\begin{enumerate}[label=\ep{\normalfont{}P\arabic*}]
		\item $f(x) = g(x)$ for all $x \in X$; and
		\item\label{item:P21} for any two colors $\alpha$ and $\beta$, if $\mu(f^{-1}(\alpha)) < \mu(f^{-1}(\beta))$, then
		\[
		\mu(\set{y \in f^{-1}(\beta) \setminus X \,:\, N_G(y) \cap f^{-1}(\alpha) = \0}) \,=\, 0.
		\]
	\end{enumerate}
	We claim that this $\Delta$-coloring $f$ is $\mu$-equitable.
	
	Suppose, towards a contradiction, that $f$ is not $\mu$-equitable. Let the set of colors be $\col{C}$. For $\gamma \in \col{C}$, let $V_\gamma \defeq f^{-1}(\gamma)$ be the corresponding color class. Let
	\[
		\col{A} \defeq \set{\alpha \in \col{C} \,:\, \mu(V_\alpha) < 1/\Delta} \qquad \text{and} \qquad \col{B} \defeq \set{\beta \in \col{C} \,:\, \mu(V_\beta) \geq 1/\Delta}.
	\]
	Set $A \defeq f^{-1}(\col{A})$ and $B \defeq f^{-1}(\col{B})$. Since $f$ is not $\mu$-equitable, $\col{A}$, $\col{B} \neq \0$. Let $\xi \defeq |\col{A}|/\Delta$. 
	
	\begin{big_claim}\label{claim:A_smallish}
		$\xi < 4/5$.
	\end{big_claim}
	\begin{scproof}
		Take any color $\beta \in \col{B}$. Using Claims~\ref{claim:X_small} and \ref{claim:h_bounds}, we see that
		\[
			\mu(V_\beta \cap X) \,=\, \mu((h^\ast)^{-1}(\beta)) \,\leq\, \frac{2\mu(X)}{\Delta+1} \,\leq\, \frac{1}{2(\Delta + 1)} \,<\, \frac{1}{\Delta} \,\leq\, \mu(V_\beta).
		\]
		Therefore, $\mu(V_\beta \setminus X) > 0$. By \ref{item:X11}, each vertex in $V_\beta\setminus X$ has degree less than $4\Delta/5$. On the other hand, by \ref{item:P21}, $\mu$-almost every vertex in $V_\beta\setminus X$ has a neighbor in every color class $V_\alpha$ with $\alpha \in \col{A}$. Thus, we have $|\col{A}| < 4\Delta/5$, or, equivalently, $\xi < 4/5$, as desired.
	\end{scproof}

	Following \cite{KN}, define $V^+ \defeq B \setminus X$ and $V^- \defeq B \cap X$.

	\begin{big_claim}\label{claim:sum}
		$\mu(B) = \mu(V^+) + \mu(V^-) \geq 1-\xi$.
	\end{big_claim}
	\begin{scproof}
		By definition, $\mu(B) \geq |\col{B}|/\Delta = 1- \xi$.
	\end{scproof}

	\begin{big_claim}\label{claim:V+}
		$\mu(V^+) < 7/10$.
	\end{big_claim}
	\begin{scproof}
		Take any $\alpha \in \col{A}$. By \ref{item:P21}, $\mu$-almost every vertex in $V^+$ has a neighbor in $V_\alpha$. Since $\mu$ is $G$-invariant, this implies that
		\[
			\int_{V_\alpha} |N_G(x) \cap V^+| \,\Diff\mu(x) \,\geq\, \mu(V^+).
		\]
		On the other hand, we can write
		\begin{align*}
			\int_{V_\alpha} |N_G(x) \cap V^+| \,\Diff\mu(x) \,&\leq\, \int_{V_\alpha \cap X} \deg_G(x) \,\Diff\mu(x) \,+\, \int_{V_\alpha \setminus X} |N_G(y) \setminus X| \,\Diff\mu(y) \\
			[\text{by \ref{item:X21}}]\qquad&\leq\, \Delta \cdot \mu(V_\alpha \cap X) \,+\, \frac{2\Delta}{5} \cdot\mu(V_\alpha \setminus X) \\
			&=\, \frac{3\Delta}{5} \cdot \mu(V_\alpha \cap X) \,+\, \frac{2\Delta}{5} \cdot \mu(V_\alpha).
		\end{align*}
		Claims~\ref{claim:X_small} and \ref{claim:h_bounds} yield $\mu(V_\alpha \cap X) \leq 1/2(\Delta + 1) < 1/2\Delta$. Also, $\mu(V_\alpha) < 1/\Delta$ since $\alpha \in \col{A}$. Thus,
		\[
			\frac{3\Delta}{5} \cdot \mu(V_\alpha \cap X) \,+\, \frac{2\Delta}{5} \cdot \mu(V_\alpha) \,<\, \frac{3\Delta}{5} \cdot \frac{1}{2\Delta} \,+\, \frac{2\Delta}{5} \cdot \frac{1}{\Delta} \,=\, \frac{3}{10} \,+\, \frac{2}{5} \,=\, \frac{7}{10},
		\]
		as desired.
	\end{scproof}

	\begin{big_claim}\label{claim:V-}
		$\mu(V^-) < (1 - \xi)/2$.
	\end{big_claim}
	\begin{scproof}
		By Claims \ref{claim:X_small} and~\ref{claim:h_bounds}, we have
		\[
			\mu(V^-) \,\leq\, \frac{(|\col{B}| + 1) \mu(X)}{\Delta+1} \,<\, \frac{|\col{B}| + 1}{4\Delta} \,\leq\, \frac{|\col{B}|}{2\Delta} \,=\, \frac{1-\xi}{2}. \qedhere
		\]
	\end{scproof}

	\begin{big_claim}\label{claim:main_ineq}
		$(4 - 10\xi) \mu(V^-) + 10\xi(1-\xi)\leq 1$.
	\end{big_claim}
	\begin{scproof}
		Observe that
		\[
			\frac{\Delta}{10} \,\geq\, \frac{d_\mu(G)}{2} \,=\, \mathsf{C}_\mu(G; V(G)) \,\geq\, \mathsf{C}_\mu(G; V^-) \,+\, \int_{V^+} |N_G(y) \cap A| \,\Diff\mu(y).
		\]
		Applying \ref{item:X31} with $X' = V^-$, we get $\mathsf{C}_\mu(G;V^-) \geq (2\Delta/5)\mu(V^-)$. Also, by \ref{item:P21}, $\mu$-almost every vertex $y \in V^+$ has at least $|\col{A}|$ neighbors in $A$. Therefore,
		\begin{equation}\label{eq:penultimate}
			\frac{\Delta}{10} \,\geq\, \frac{2\Delta}{5} \cdot \mu(V^-) \,+\, |\col{A}| \cdot \mu(V^+).
		\end{equation}
		By Claim~\ref{claim:sum},
		\[
			\frac{2\Delta}{5} \cdot \mu(V^-) \,+\, |\col{A}| \cdot \mu(V^+) \,\geq\, \left(\frac{2\Delta}{5} - |\col{A}|\right) \cdot \mu(V^-) \,+\, |\col{A}| \cdot (1-\xi).
		\]
		Plugging this into \eqref{eq:penultimate} and multiplying both sides by $10/\Delta$ gives the desired result.
	\end{scproof}

	\begin{big_claim}
		$\xi < 2/5$.
	\end{big_claim}
	\begin{scproof}
		Suppose that $\xi \geq 2/5$. Then $4-10\xi \leq 0$, so it follows from Claim~\ref{claim:main_ineq} that
		\begin{align*}
			1 \,&\geq\, (4 - 10\xi) \mu(V^-) \,+\, 10\xi(1-\xi) \\
			[\text{by Claim~\ref{claim:V-}}]\qquad&\geq\, (2-5\xi)(1-\xi) \,+\, 10\xi(1-\xi).
		\end{align*}
		In other words, we have $5\xi^2 - 3\xi - 1 \geq 0$. This inequality implies that either $\xi \leq (3-\sqrt{29})/10 < 0$, which is impossible, or else, $\xi \geq (3+\sqrt{29})/10 = 0.83\ldots > 4/5$, contradicting Claim~\ref{claim:A_smallish}.
	\end{scproof}

	We are ready for the coup de gr\^{a}ce. Since $\xi < 2/5$, we have $4-10\xi > 0$, so Claim~\ref{claim:main_ineq} yields
	\[
		\mu(V^-) \,\leq\, \frac{1 - 10\xi(1-\xi)}{4-10\xi} \,=\, \frac{10\xi^2 - 10\xi + 1}{4-10\xi}.
	\]
	From this and Claim~\ref{claim:sum} we obtain
	\[
		\mu(V^+) \,\geq\, 1 - \xi - \frac{10\xi^2 - 10\xi + 1}{4-10\xi} \,=\, \frac{3-4\xi}{4-10\xi}.
	\]
	The function $t \mapsto (3-4t)/(4 - 10t)$ is increasing for $0 \leq t < 2/5$, so $\mu(V^+)$ is at least the value of this function at $t = 0$, i.e., $\mu(V^+) \geq 3/4$. Since $3/4 > 7/10$, this contradicts Claim~\ref{claim:V+} and completes the proof of Theorem~\ref{theo:meas_KN}.

	{\renewcommand{\markboth}[2]{}
		\printbibliography}

\end{document}